\documentclass[12pt]{amsart}
\usepackage{amsfonts}
\usepackage{amsmath}
\usepackage{amssymb}
\usepackage[margin=1in]{geometry}
\usepackage{tikz}
\usepackage{graphicx}
\usepackage{amsthm}   							% hiermee kan je zelf omgevingen zoals Stelling definiëren
\usepackage{mathtools}							% uitbreiding van amsmath 
\usepackage{commath}							% Om differentialen te typesetten als operatoren:
\usepackage{enumitem}  							% hiermee kan je zelf de labels kiezen bij vb. itemize en enumeratie

\usepackage[english]{babel}						% http://ftp.snt.utwente.nl/pub/software/tex/macros/latex/required/babel/base/babel.pdf

\makeatletter
\@namedef{subjclassname@2020}{\textup{2020} Mathematics Subject Classification}
\makeatother

\newtheorem{theorem}{Theorem}[section]

\newtheorem{lemma}[theorem]{Lemma}
\newtheorem{proposition}[theorem]{Proposition}

\theoremstyle{definition}

\newtheorem{remark}[theorem]{Remark}

\numberwithin{equation}{section}

\newcommand{\N}{\mathbb{N}}
\newcommand{\Z}{\mathbb{Z}}

\newcommand{\R}{\mathbb{R}}
\newcommand{\C}{\mathbb{C}}

\renewcommand{\Re}{\operatorname{Re}}
\renewcommand{\Im}{\operatorname{Im}}

\newcommand{\I}{\mathrm{i}}
\newcommand{\e}{\mathrm{e}}

\newcommand{\eps}{\varepsilon}
\newcommand{\vphi}{\varphi}

\newcommand{\D}{\mathcal{D}}
\newcommand{\MS}{\mathcal{S}}
\newcommand{\ML}{\mathcal{L}}
\newcommand{\F}{\mathcal{F}}
\newcommand{\K}{\mathcal{K}}

\newcommand{\fracD}[3]{\prescript{}{#1}{D}^{#2}_{#3}}	% Fractional derivative with initial point #1, order of differentiation #2, and variable #3.

\DeclareMathOperator{\supp}{supp}
\DeclareMathOperator{\singsupp}{sing\, supp}
\DeclareMathOperator{\WF}{WF}
\DeclareMathOperator{\sgn}{sgn}
\DeclareMathOperator{\pv}{p.v.}
\DeclareMathOperator{\Fp}{F.p.}
\DeclareMathOperator{\length}{length}
\DeclareMathOperator{\argmax}{argmax}

\begin{document}

\title{Micro-local and qualitative analysis of the fractional Zener wave equation}

\author[F. Broucke]{Frederik Broucke}
\thanks{Frederik Broucke was supported by the Ghent University BOF-grant 01J04017}
\address{Department of Mathematics: Analysis, Logic and Discrete Mathematics\\ Ghent University\\ Krijgslaan 281\\ 9000 Gent\\ Belgium}
\email{fabrouck.broucke@ugent.be, ORCID: 0000-0002-5744-4767}

\author[Lj. Oparnica]{Ljubica Oparnica}
\thanks{Ljubica Oparnica was supported by the FWO Odysseus 1 grant no. G.0H94.18N}
\address{Department of Mathematics: Analysis, Logic and Discrete Mathematics, Ghent University, Krijgslaan 281,  9000 Gent, Belgium}

\email{oparnica.ljubica@ugent.be,  ORCID: 0000-0001-8547-339X}

\begin{abstract}
This paper concerns the micro-local and qualitative analysis of the fractional Zener wave equation.
The classical and Gevrey-type wave front sets of the fundamental solution are determined, and  questions on dispersion, dissipation, wave propagation speed and wave packet shape are addressed. 
\end{abstract}

%35B65 Smoothness and regularity of solutions to PDEs
%35R11 Fractional partial differential equations
%74J05 Linear waves in solid mechanics
%74D05 Linear constitutive equations for materials with memory
%41A60 Asymptotic approximations, asymptotic expansions (steepest descent, etc.)

\maketitle

%\tableofcontents
%%%%%%%%%%%%%%%%%%%%%%%%%%%%%%%%%%%%%%%%%%%%%%%%%%%%%%%%%%%%%%%%%%%%%%%%%%%%%
% INTRODUCTION %%%%%%%%%%%%%%%%%%%%%%%%%%%%%%%%%%%%%%%%%%%%%%%%%%%%%%%%%%%%%%%%%%%
%%%%%%%%%%%%%%%%%%%%%%%%%%%%%%%%%%%%%%%%%%%%%%%%%%%%%%%%%%%%%%%%%%%%%%%%%%%%%
\section{Introduction}

In this paper, we study the fractional Zener wave equation, which models the propagation of waves in a viscoelastic medium and which is given by
\begin{equation}
\label{eq: FZWE}
	\dpd[2]{}{t} u(x,t) = \ML^{-1}_{s\to t} \left(\frac{1+s^{\alpha}}{1+\tau s^{\alpha}}\right) \ast_{t} \dpd[2]{}{x} u(x,t) ,\qquad x\in \mathbb{R}, \quad t>0,
\end{equation}
where $\ML^{-1}$ denotes the inverse Laplace transform, while $\tau$ and $\alpha$ are constants satisfying $0<\tau<1$, $0<\alpha<1$.

%Wave propagation in a viscoelastic medium is modelled via the system of basic equations of continuum mechanics, in which the equation of motion and the equation for the strain describing small local deformations are unchanged, compared to the purely elastic case, while Hooke's law, the constitutive law for an elastic body, is replaced by a law which describes the mechanical properties of the linear viscoelastic body, see \cite{a-g, APSZ-2}. 
Like the classical wave equation, the fractional Zener wave equation can be derived from a system of three basic equations of continuum mechanics. The equation of motion and the equation for the strain describing small local deformations are unchanged compared to the purely elastic case, while Hooke's law, the constitutive law for an elastic body, is replaced by a law which describes the mechanical properties of the linear \emph{viscoelastic} body. For more details, we refer to \cite{a-g, APSZ-2, Mainardi}. The constitutive law yielding the fractional Zener wave equation \eqref{eq: FZWE} is the fractional Zener law, proposed in papers of Caputo and Mainardi \cite{CaputoMainardi-1971b,CaputoMainardi-1971a}. Some restrictions on the coefficients appearing in this constitutive equation, necessary for being physically consistent with the second law of thermodynamics and known as thermodynamical restrictions, are given in \cite{BT1986}. 
In \cite{KOZ2010}, the system with the fractional Zener law instead of the Hooke's law was shown to be equivalent to the fractional Zener wave equation \eqref{eq: FZWE}. The existence and uniqueness of solutions to the generalized Cauchy problem in the space of tempered distributions was proven under the assumption of the thermodynamical restrictions, a representation of the solution was given, and numerical examples were provided indicating some qualitative properties of the solutions, such as properties concerning shape and dissipation.
%In \cite{APS}, the fractional Zener wave equation was also investigated, and a statistical analysis was performed for solutions with stochastic Cauchy data.
In \cite{ZO2020}, energy dissipation was proven for a general class of fractional wave equations which includes the fractional Zener wave equation.
The micro-local analysis for the fractional wave equation was initiated in \cite{HOZ16}, where an analogue of non-characteristic regularity was shown.
% but the wave front set for the fractional Zener wave equation was not precisely determined.

The work we are to present continues the analysis of the fractional Zener wave equation, 
%Should we mention the support in the intro?
%proving rigorously the support properties of the fundamental solution already announced earlier in  \cite{KOZ2010}, 
completing and extending the results on the wave front set given in \cite{HOZ16}. We provide a complete description of the $C^{\infty}$-wave front set of the fundamental solution of \eqref{eq: FZWE}. In particular we show that the fundamental solution is smooth on the boundary of the forward light cone, in contrast with the classical wave equation. We also determine the wave front set with respect to Gevrey classes of functions, which assumes a finer notion of smoothness. As a consequence we prove that for the order of the Gevrey class sufficiently close to 1, the fundamental solution is singular at the boundary of the forward light cone.

Next, we perform a qualitative analysis of solutions in two cases. First, we investigate the response of the system when it is submitted to a forced harmonic oscillation at the origin. From this we detect the presence of dissipation and anomalous dispersion. Second, we investigate the evolution following a ``delta impulse'', i.e.\ the solution with initial conditions $u(x,0)=\delta(x)$, $\partial_{t}u(x,0) = 0$ (and zero force term). This solution consists of two wave packets moving in opposite directions. We will provide an accurate description of the ``limiting shape'' of this wave packet, and motivate the notion of a wave packet speed. 
%and providing a qualitative analysis with implications to questions on dispersion and energy dissipation, wave speed and wave packet shape, being of particular importance for applications in viscoelasticity. The content of the paper is organized as follows.

The content of the paper is organized as follows.
The first section is introductory and provides motivation and mathematical preliminaries necessary for the work. In order to make the paper self-contained, the set up of the Cauchy problem for the fractional Zener wave equation under consideration with some details from previous works is given at the beginning of Section \ref{FS}. We further discuss representation formulas and properties of the fundamental solution. 
% where afterwards, the fundamental solution of the fractional Zener wave equation is analysed  with the main result given in Proposition \ref{prop: support S} proving that the fundamental solution is supported in a forward cone.  
The main result of the paper concerning the regularity of the fundamental solution, is given in Section \ref{MLA}. Theorem \ref{thm: regularity S} and Theorem \ref{th: Gevrey regularity} describe the wave front set of the fundamental solution with respect to $C^{\infty}$ and the Gevrey classes respectively.
%proves that in non-characteristic regularity inclusion relation given in \cite[Theorem 3.2]{HOZ16}, the points in the boundary of the forward light cone in fact do not belong to wave front set and that the fundamental solution of fractional Zener wave equation is smooth on the boundary of the forward light cone.
 %Since this is in contrast with the classical wave equation regularity property, we investigate singularities by the use of Gevrey classes of functions, which assumes a finer notion of smoothness. The result given in Theorem \ref{th: Gevrey regularity}  shows that for the order of the Gevrey class sufficiently close to 1, the boundary of the forward light cone becomes singular. 
The proof of the latter theorem is long and quite technical. The main ideas are presented there, but the proof of a technical lemma is provided in Appendix \ref{sec: technical calculation}. Section \ref{QA} concerns the qualitative analysis and is divided into two subsection, treating the forced harmonic oscillation and the delta impulse respectively.
%In the first subsection, the response when the system is submitted to a forced oscillation is investigated, and the presence of dissipation and anomalous dispersion is shown. In the second subsection, the shape of solutions with Cauchy data $u(x,0)=u_{0}(x)$, $\partial_{t}u(x,0) = 0$ and $f(x,t) =0$ is analysed and described. 
Section \ref{ClassZ} addresses the case when the viscoelastic medium is described by the \emph{classical} Zener model (or the Standard Linear Solid (SLS) model) and shows that, as for the classical wave equation and in contrast to the fractional Zener wave equation, the fundamental solution is not smooth on the boundary of the light cone. Moreover, some qualitative difference between the two models in terms of dissipation is noted.
% in the SLS model, two pseudo-monochromatic waves with different frequencies have roughly the same amount of spatial dampening, while in the fractional Zener model, the wave with the higher frequency will experience more dampening than the wave with the lower frequency.

%%%%%%%%%%%%%%%%%%%%%%%%%%%%%%%%%%%%%%%%%%%%%%%%%%%%%%%%%%%%%%%%%%%%%%%%%%%%%
% MATHEMATICAL PRELIMINARIES %%%%%%%%%%%%%%%%%%%%%%%%%%%%%%%%%%%%%%%%%%%%%%%%%%%%%%%%%
%%%%%%%%%%%%%%%%%%%%%%%%%%%%%%%%%%%%%%%%%%%%%%%%%%%%%%%%%%%%%%%%%%%%%%%%%%%%%

\subsection{Preliminaries}

\subsubsection{Notations}
 For  $\Omega$  an open subset of $\R^{n}$, $\mathcal{D}(\Omega)$ denotes the space of compactly supported smooth functions on $\Omega$. The space of distributions on $\Omega$ is denoted by $\mathcal{D}'(\Omega)$ and the space of compactly supported distributions is denoted by $\mathcal{E}'(\Omega)$. We denote with $\mathcal{S}(\R^{n})$ the Schwartz space of rapidly decreasing smooth functions, and with $\mathcal{S}'(\R^{n})$ the space of tempered distributions. Further, $\mathcal{D}'_{+}(\R) \subset \mathcal{D}'(\R)$ and $\mathcal{S}'_{+}(\R) \subset \mathcal{S}'(\R)$ are spaces of distributions supported on $[0,\infty)$, and $\mathcal{S}'(\mathbb{R}\times \mathbb{R}_{+})$ is the space of  distributions in $\mathcal{S}'(\mathbb{R}^{2})$ vanishing on $\R\times(-\infty, 0)$.
  
We use the following common notations for asymptotic relations:
\begin{align*}
	f(x) &\lesssim g(x) \iff f(x) = O(g(x)) \iff \exists C>0\colon \abs{f(x)} \le Cg(x); \\
	f(x) &\simeq g(x) \iff f(x)\lesssim g(x) \text{ and } g(x) \lesssim f(x); \\
	f(x) &\sim g(x) \iff \frac{f(x)}{g(x)}\to 1.
\end{align*}
The denote the dependency of implicit $O$-constants on parameters, we will use subscripts, for example $\lesssim_{n}$.
 
% Since we shall mostly work with functions and distributions depending on two variables, $u=u(x,t)$, we introduce 
% We shall also need the set $\mathbb{C}_{+}=\{z\in\mathbb{C}\,|\, \Re z>0\}$.

\subsubsection{Integral transforms}We define the Fourier transform for an integrable function $\varphi \in L^1(\R)$
%Schwartz function  $\varphi \in \mathcal{S}(\mathbb{R})$ 
 as
\[
	 \F\varphi (\xi)=\hat{\varphi}(\xi)=\int_{-\infty }^{\infty}\varphi (x)\e^{-\I\xi x}\dif x, \qquad \xi \in \mathbb{R},
\]
while for  $f\in \mathcal{S}'(\mathbb{R})$ the Fourier transform is given via $\left\langle \F f,\varphi \right\rangle =\left\langle f,\F\varphi \right\rangle$, $\varphi \in \mathcal{S}(\mathbb{R})$.

 The Laplace transform of $f\in \mathcal{D}'_{+}(\mathbb{R})$ satisfying $\e^{-at}f\in \mathcal{S}'(\mathbb{R})$, for all $a >a_{0}>0$ is defined by
\[
 \ML f(s) =\widetilde{f}(s) = \F(\e^{-a t}f)(y), \qquad s=a+\I y,
\]
$\ML f$ being a holomorphic function in the half plane $\Re s>a_{0}$. %(see e.g. \cite{DautryLions-vol5} or \cite{Vladimirov-EqMP}). 
If $f\in\MS'(\R_{+})$, then $\ML f(s) = \langle f(t), \e^{-st}\rangle$. In particular, for 
$f\in L^1(\mathbb{R})$ with $f(t)=0$, for $t<0$, %and $|y(t)| \leq Ae^{at}$ ($a,A>0$), 
the Laplace transform is given by
\[
 \ML f(s)=\int_{0}^{\infty} f(t)\e^{-st}\dif t, \qquad \Re s\ge0.
\]
The inverse Laplace transform exists as distribution in $\MS_{+}(\R)$ for functions $F$  holomorphic in the half  plane $\Re s>a_{0}$ satisfying $\abs{F(s)} \leq A\frac{(1+\,\abs{s})^{m}}{\abs{\Re s}^{k}}$, $m,k\in \N$, $\Re s>0$, and is given by

\begin{equation} 
\label{eq: Laplace inverse formula}
 f(t) = \ML^{-1}F(t) =\lim_{Y\to\infty}\frac{1}{2\pi\I}\int_{a-\I Y}^{a+\I Y} F(s) \e^{st}\dif s, \qquad t>0, \quad a>a_{0},
\end{equation}
whenever this limit exists.

If $f(x,t) \in \MS'(\R\times\R_{+})$, then the Laplace transform of $f$ with respect to $t$ is the distribution-valued function
\[
	\ML_{t}f: \{s: \Re s>0\} \to \MS'(\R): s \mapsto \bigl(\phi(x) \mapsto \langle f(x,t), \phi(x)\e^{-st}\rangle\bigr).
\]
\subsubsection{Fractional derivatives}
The left  Riemann-Liouville fractional derivative of order $\alpha\in [0,1)$ is defined for an absolutely continuous function, $f\in AC([0,a])$, on an interval $[0,a]$ with $a>0$ by
\[
 	\fracD{0}{\alpha}{t}f(t) =\frac{1}{\Gamma (1-\alpha)} \od{}{t} \int_{0}^{t} \frac{f(\zeta)}{(t-\zeta)^{\alpha}}\dif\zeta, \quad t\in [0,a],
\]
 while the left Liouville-Weyl fractional derivative of order $\alpha\in [0,1)$ is defined for $f\in AC(\mathbb{R})$ with $f(-t) \lesssim 1/t$ for $t\to \infty$ by
\[
	\fracD{-\infty}{\alpha}{t} f(t) =\frac{1}{\Gamma (1-\alpha)} \od{}{t} \int_{-\infty}^{t} \frac{f(\zeta)}{(t-\zeta)^{\alpha}}\dif\zeta, \quad t\in \mathbb{R},
\]
 where $\Gamma$ denotes the Euler gamma function. In the spaces of distributions one introduces a family $\{\chi_{+}^{\alpha}\}_{\alpha\in\C}\in \MS'_{+}(\R)$:
\[
	\chi_{+}^{\alpha}(t) = \frac{1}{\Gamma(\alpha+1)}t_{+}^{\alpha}.
\]
This family of distributions is initially defined for $\Re\alpha>-1$ as $L^{1}_{\mathrm{loc}}$ functions, but can be extended to every $\alpha\in\C$ by analytic continuation. When $\alpha$ is not a negative integer, $\langle\chi_{+}^{\alpha}(t), \vphi(t)\rangle$ can be evaluated using the Hadamard finite part:
\[
	\langle\chi_{+}^{\alpha}(t), \vphi(t)\rangle = \frac{1}{\Gamma(\alpha+1)}\Fp\int_{0}^{\infty}t^{\alpha}\vphi(t)\dif t, \quad \vphi\in \MS(\R).
\]
For more details concerning this family of distributions, we refer to \cite[Section 3.2]{Hormander}.

The convolution operator $f\mapsto \chi_{+}^{-\alpha-1}\ast f$ coincides with the left Riemann-Liouville and left Liouville-Weyl fractional derivative of order $\alpha$, for $f\in AC([0,a])$ and $f\in AC(\R)$ respectively. 
 
The Fourier transform of the Liouville-Weyl fractional derivative of $f\in AC(\mathbb{R})$ and Laplace transform of the left Riemann-Liouville fractional derivative of $f\in AC([0,a])$ are given as follows:
\[
	 \F[\fracD{-\infty}{\alpha}{t} f](\xi) =( \I\xi)^{\alpha} \F f(\xi),  \qquad \ML[\fracD{0}{\alpha}{t}f](s) = s^{\alpha}\ML f(s).
\]
For more details on fractional derivatives we refer to \cite{SKM}.

\subsubsection{Gevrey classes} Let $\sigma\ge 0$ and let $\Omega \subseteq \R^{n}$ be open. A function $\vphi\in C^{\infty}(\Omega)$ belongs to the Gevrey class $G^{\sigma}(\Omega)$ of order $\sigma$ if for every compact $K\subseteq\Omega$ there exists a constant $C=C_{K}>0$ such that 
\[
	\sup_{x\in K}\abs[1]{\partial^{\beta}\vphi(x)} \le C^{1+\, \abs{\beta}}(\beta!)^{\sigma}, \quad \text{for every multi-index } \beta\in\N^{n}.
\]
The case $\sigma=1$ corresponds to real analytic functions. When $\sigma>1$, $G^{\sigma}(\Omega)$ contains compactly supported functions. For a distribution $u\in\D'(\Omega)$, $\singsupp_{G^{\sigma}} u$ is defined as the complement of the largest open set $X$ where $u\in G^{\sigma}(X)$. Similarly as in the $C^{\infty}$ case, one might perform a spectral analysis of the $G^{\sigma}$-singularities of $u$, by investigating the decay on cones of the Fourier transform of localizations of $u$. Since the space $G^{\sigma}$, $\sigma\le1$ does not contain compactly supported functions, one uses so-called analytic cut-off sequences to localize. This leads to the notion of the $G^{\sigma}$-wave front set of $u$, denoted as $\WF_{G^{\sigma}}(u)$. For more details we refer to \cite[Section 8.4]{Hormander}.

\subsubsection{The saddle point method} To estimate certain integrals we will use the saddle point method, also known as the method of steepest descent. This method can be summarized as follows. Suppose $f$ and $g$ are holomorphic functions on a simply connected region $\Omega$. Suppose that $a, b \in \Omega$, and suppose that $z_{0}$ is a simple saddle point of $f$, that is, a simple zero of $f'$. If there exists a path from $a$ to $b$ in $\Omega$, which passes through $z_{0}$ in such a way that $\Re f(z)$ reaches its maximum at $z=z_{0}$, then 
\[
	\int_{a}^{b}g(z)\e^{\lambda f(z)}\dif z \sim g(z_{0})\sqrt{\frac{2\pi}{-f''(z_{0})\lambda}}\e^{\lambda f(z_{0})}, \quad \text{as } \lambda \to \infty.
\]
For an introduction to the saddle point method, we refer to \cite[Chapter 5]{deBruijn}.

%\subsubsection{}\todo{Maybe something...}

\section{The Cauchy problem and the fundamental solution}\label{FS}
%%%%%%%%%%%%%%%%%%%%%%%%%%%%%%%%%%%%%%%%%%%%%%%%%%%%%%%%%%%%%%%%%%%%%%%%%%%%%
% THE CAUCHY PROBLEM %%%%%%%%%%%%%%%%%%%%%%%%%%%%%%%%%%%%%%%%%%%%%%%%%%%%%%%%%%%%%%
%%%%%%%%%%%%%%%%%%%%%%%%%%%%%%%%%%%%%%%%%%%%%%%%%%%%%%%%%%%%%%%%%%%%%%%%%%%%%

To describe waves occurring in one-dimensional viscoelastic media one uses the system of basic equations of elasticity (see \cite{a-g}), consisting of the equilibrium equation coming from Newton's second law: ${\partial_x}\sigma = \rho \,{\partial^{2}_t} u$, the constitutive equation, describing the connection between stress and strain, and the strain measure for local small deformations, giving the connection between strain and displacement: $\varepsilon =\partial_x u$. Here, $\sigma$, $u$ and $\varepsilon$ denote stress, displacement, and strain respectively, being functions of $x\in \mathbb{R}$ and $t>0$. The constant $\rho$ denotes the density of the medium under consideration. 

For purely elastic media, the constitutive equation is given by Hooke's law: $\sigma = E\eps$, where $E$ is a constant referred to as the Young modulus of elasticity. The classical wave equation $\partial^{2}_t u  = c^{2}\partial^{2}_x u$  is then easily obtained from this system with $c=\sqrt{E/\rho}$ being a constant that can be physically interpreted as the wave speed.

For viscoelastic media one replaces Hooke's law by other constitutive equations, such as the Maxwell, Voight, or Zener models (see \cite{Mainardi}).
 In the fractional models, one replaces the integer order derivatives by derivates of  fractional order. 
 The fractional Zener model is a generalization of the (classical) Zener model, also known as the Standard Linear Solid (SLS) model, and is given by
% The fractional Zener constitutive equation is given by
 \begin{equation} 
 \label{eq:frac Zener model}
 	\sigma (x,t) +\tau_{\sigma} \fracD{}{\alpha}{}\sigma (x,t) = E[\varepsilon (x,t) +\tau_{\varepsilon} \fracD{}{\alpha}{}\varepsilon (x,t)],
 \qquad x\in \mathbf{\mathbb{R}},\quad t>0, 
 \end{equation}
 where $\fracD{}{\alpha}{}$ denotes fractional differentiation of order $\alpha$, $0<\alpha<1$. The case $\alpha=1$, i.e. with ordinary derivatives, corresponds to the SLS model.
 In \eqref{eq:frac Zener model}, $\tau_{\sigma}$, $\tau_{\varepsilon}$ are constants referred to as the relaxation time and the retardation time respectively. They satisfy the thermodynamical restriction $\tau_{\varepsilon}>\tau_{\sigma}>0$ following from the second law of thermodynamics \cite{BT1986}. The constant $E$ is the Young modulus of elasticity.

In this paper, we analyze wave propagation in viscoelastic media modelled by the system of basic equations  in which the fractional Zener law \eqref{eq:frac Zener model} instead of Hooke's law is used. In dimensionless form, this system is given by (see e.g. \cite{KOZ2010})
%We shall work with the system reduced to dimensionless quantities, given by (see \cite{KOZ2010} for the details)
\begin{align}
	\dpd{}{x}\sigma 						&= \dpd[2]{}{t} u,  \nonumber \\
 	\sigma + \tau \fracD{0}{\alpha}{t}\sigma 	&= \varepsilon +\fracD{0}{\alpha}{t}\varepsilon,  \label{eq:system} \\
	\varepsilon 						&=  \dpd{}{x} u; \nonumber
\end{align}
 where as above, $\sigma$, $u$ and $\varepsilon$ are stress, displacement, and strain, respectively, considered as functions of $x\in \mathbb{R}$ and $t>0$, $0<\tau<1$ is
 a constant and $\fracD{0}{\alpha}{t}$, $0< \alpha <1$, is the left Riemann-Liouville operator of fractional differentiation. In \cite{KOZ2010} it was shown that the wave equation corresponding to \eqref{eq:system} is the fractional Zener wave equation \eqref{eq: FZWE}.
% as well as the boundary conditions
%\[
%	 \lim_{x\to\pm \infty}u(x,t)=0, \qquad \lim_{x\to \pm \infty} \sigma(x,t)=0.
%\]
%\todo{What if $u_{0}\in \MS'$? Then not necessarily $\lim u = 0$?}
%In \cite{KOZ2010}, it was shown that the system \eqref{eq:system} can be reduced to fractional Zener wave equation \eqref{eq: FZWE}.
% \begin{equation*} %\label{eq:fzwe}
%\dpd[2]{}{t} u(x,t) = \ML^{-1}_{s\to t} \left(\frac{1+s^{\alpha}}{1+\tau s^{\alpha}}\right) \ast_{t} \dpd[2]{}{x} u(x,t) ,\qquad x\in \mathbb{R}, t>0.
% \end{equation*}
The convolution kernel $\ML^{-1}\bigl((1+s^{\alpha})/(1+\tau s^{\alpha})\bigr)$ can be expressed using Mittag-Leffler functions (see \cite[Appendix E]{Mainardi}) as
\[
	\ML^{-1}\biggl(\frac{1+s^{\alpha}}{1+\tau s^{\alpha}}\biggr) = \frac{1}{\tau}\delta(t) - \frac{1-\tau}{\tau^{2}}e_{\alpha,\alpha}(t; 1/\tau)= \frac{1}{\tau}\delta(t) + \frac{1-\tau}{\tau} \dod{}{t} e_{\alpha}(t; 1/\tau).
\] 
% {\partial_x}\sigma = \rho \,{\partial^{2}_t} u
%\todo{connect}
Further in \cite{KOZ2010}, 
%the system \eqref{eq:system} was subjected to initial conditions \todo{But $\eps = \partial_{x}u = u_{0}'(x)$?}
%\[  
%	u(x,0) = u_{0}(x), \quad \pd{u}{t}(x,0) = v_{0}(x), \quad \sigma (x,0) = 0, \quad \varepsilon(x,0) = 0,
%\]
%and 
for given initial data  $u_{0}$, $v_{0} \in \MS'(\R)$, and force term $f\in\MS'(\R\times\R_{+})$, the generalized Cauchy problem for the fractional Zener wave equation% given by 
\begin{equation}
\label{eq: Cauchy problem Zener}
Pu(x,t) = f(x,t) + u_{0}(x)\delta'(t) + v_{0}(x)\delta(t),
\end{equation}
with
\[ P = \pd[2]{}{t} - \ML^{-1}\biggl(\frac{1+s^{\alpha}}{1+\tau s^{\alpha}}\biggr) \ast_{t} \pd[2]{}{x}, \] 
was considered. It was shown that \eqref{eq: Cauchy problem Zener} has a unique solution expressed via convolution of the fundamental solution $S$ with the Cauchy data:
\begin{equation}
\label{eq: solution}
	u(x,t) = S(x,t)\ast\bigl( f(x,t) + u_{0}(x)\delta'(t) + v_{0}(x)\delta(t)\bigr).
\end{equation}
The Laplace transform $\tilde{S}$ of the fundamental solution $S$ with respect to $t$ is calculated as
\begin{equation}
\label{eq: tilde{S}}
	\tilde{S}(x,s) = \frac{1}{2s}\sqrt{\frac{1+\tau s^{\alpha}}{1+s^{\alpha}}}\exp\biggl(-\abs{x}s\sqrt{\frac{1+\tau s^{\alpha}}{1+s^{\alpha}}}\biggr), \quad x\in \R, \, \Re s > 0,
\end{equation}
where the principal branch of the logarithm is used for the function $s^{\alpha}$ and the square root. Note that for fixed $s$, this is a continuous function of $x$. Denote by $l_{\alpha}(s)$ the function defined as 
\begin{equation*}
	l_{\alpha}(s) = \sqrt{\frac{1+\tau s^{\alpha}}{1+s^{\alpha}}}, \quad \arg s \in [-\pi,\pi].
\end{equation*}
First we derive some properties of  $l_{\alpha}(s)$.
\begin{lemma}
\label{lem: l_{alpha}}
The real and imaginary part of $l_{\alpha}$ satisfy
\begin{equation}
\label{eq: sgn l_{alpha}}
	\Re l_{\alpha}(s) > 0, \quad \sgn \Im l_{\alpha}(s) = -\sgn \Im s, \quad \arg s \in [-\pi,\pi]. 
\end{equation}
Its asymptotic behavior near the origin and infinity is given by
\begin{align}
	l_{\alpha}(s)	&= 1- \frac{1-\tau}{2}s^{\alpha} + O(\,\abs{s}^{2\alpha}),  									&&\text{as } \abs{s} \to 0; \label{eq: approx l_{alpha} 0}\\
	l_{\alpha}(s)	&= \sqrt{\tau}\biggl(1 + \frac{1}{2}\biggl(\frac{1}{\tau}-1\biggr)s^{-\alpha} + O(\,\abs{s}^{-2\alpha})\biggr),  	&&\text{as } \abs{s}\to\infty. \label{eq: approx l_{alpha}}
\end{align}
\end{lemma}
In particular, there exist positive constants $c_{1}$ and $c_{2}$ such that
\begin{equation}
\label{eq: bound l_{alpha}}
	\Im l_{\alpha}(\I y) \le \begin{dcases}
		-c_{1}y^{\alpha}		&\text{for } 0\le y \le 1;\\
		-c_{2}y^{-\alpha}	&\text{for } y\ge1.
	\end{dcases}
\end{equation}
\begin{proof}
A straightforward calculation shows that, with $s=R\e^{\I\vphi}$, $-\pi\le\vphi\le\pi$,
\[
	\frac{1+\tau s^{\alpha}}{1+s^{\alpha}} = \frac{1+\tau R^{2\alpha} + (1+\tau) R^{\alpha}\cos(\alpha\vphi) 
	-\I (1-\tau) R^{\alpha}\sin(\alpha\vphi)}{1+R^{2\alpha}+2R^{\alpha}\cos(\alpha\vphi)}.
	%R^{\alpha}(1-\tau) \e^{-\I\alpha\vphi}}.
\]
The denominator is real and positive. We see that $\frac{1+\tau s^{\alpha}}{1+s^{\alpha}} \in \C \setminus (-\infty, 0]$, so the real part of its square root is positive. Since taking the square root does not alter the sign of the imaginary part, the first claim of the lemma follows. The formulas \eqref{eq: approx l_{alpha} 0} and \eqref{eq: approx l_{alpha}} follow immediately from Taylor's formula, and  upon writing $l_{\alpha}(s) = \sqrt{\tau}\sqrt{(1+\tau^{-1}s^{-\alpha})/(1+s^{-\alpha}})$ for large $s$.
\end{proof}

In \cite{KOZ2010}, the following representation of $S$ inside the forward cone $\abs{x}<t/\sqrt{\tau}$ was given\footnote{The constant right after the equality sign in \cite[eq.\ (18)]{KOZ2010} should be $1/2$ instead of $1$.}
\begin{equation}
\label{eq: S in cone}
	S(x,t) = \frac{1}{2} + 
	\frac{1}{4\pi\I}\int_{0}^{\infty}\bigl(l_{\alpha}(q\e^{\I\pi})\e^{\,\abs{x}ql_{\alpha}(q\e^{\I\pi})} 
	- l_{\alpha}(q\e^{-\I\pi})\e^{\,\abs{x}ql_{\alpha}(q\e^{-\I\pi})}\bigr)\frac{\e^{-qt}}{q}\dif q.
\end{equation}
Note that it follows from the asymptotic behavior of $l_{\alpha}$ that this integral converges absolutely whenever $\abs{x} < t/\sqrt{\tau}$.
The representation \eqref{eq: S in cone} was shown by Laplace inversion via formula \eqref{eq: Laplace inverse formula}, i.e.\ by calculation of the integral
\begin{equation} 
\label{eq: Laplace inverse}
\frac{1}{2\pi\I}\lim_{Y\to\infty}\int_{a-\I Y}^{a+\I Y}\frac{ l_{\alpha}(s)}{2s}  \e^{-\, \abs{x} s l_{\alpha}(s)+ts} \dif s, \quad a>0,
\end{equation}

%For $\abs{x}<t/\sqrt{\tau}$, it was shown that the above integral converges by switching the contour of integration to a Hankel contour encircling the negative real axis. The contribution of the singularity at $s=0$ is $1/2$, and combining the integrals along both sides of the branch cut $(-\infty, 0]$ gives the absolutely convergent integral in \eqref{eq: S in cone}. 

%If $\abs{x}>t/\sqrt{\tau}$, it turns out that $S(x,t)=0$. 
There is some discrepancy in the literature (e.g.\   \cite{APS} and \cite{KOZ2010}) regarding claims about the support of $S$. We will clarify this here and prove\footnote{We note that support properties for a more general class of models were also proven in \cite{KOZ2019}} that $S$ is supported in the forward cone $\abs{x}\le t/\sqrt{\tau}$. We will also indicate how to deduce the representation \eqref{eq: S in cone}, since this technique will be used multiple times throughout this paper.
%This was claimed in \cite{KOZ2010} without proof while in \cite{KOZ2019}, a proof is provided for a more general class of models. 
%For convenience of the reader, we will provide a proof here. 
\begin{proposition}
\label{prop: support S}
The fundamental solution $S$ of \eqref{eq: FZWE} is supported in a forward cone:
\[
	\supp S \subseteq \{(x,t): \abs{x} \le t/\sqrt{\tau}\}.
\] 
In the interior of this cone, $S$ is given by \eqref{eq: S in cone}.
\end{proposition}
\begin{proof}%[Proof of Proposition \ref{prop: support S}]
Let $x$ and $t$ be such that $\abs{x}>t/\sqrt{\tau}$. We show that $S(x,t)=0$. Using Cauchy's formula, we may rewrite the integral \eqref{eq: Laplace inverse} as an integral over the arc of the circle of radius $R=\sqrt{a^{2}+Y^{2}}$ and center $0$, which connects the points $a-\I Y$ and $a+\I Y$. The  polar angle varies between $-\vphi(R)$ and $\vphi(R)$, with $\vphi(R)=\arctan\sqrt{(R/a)^{2}-1}$. We get
\[
	S(x,t) = \lim_{R\to\infty}\frac{1}{4\pi}\int_{-\vphi(R)}^{\vphi(R)}l_{\alpha}(R\e^{\I\vphi})\exp\bigl(R\e^{\I\vphi}(t-\abs{x}l_{\alpha}(R\e^{\I\vphi}))\bigr)\dif\vphi.
\]
Using \eqref{eq: sgn l_{alpha}} and extending the range of integration to $[-\pi/2,\pi/2]$, we can bound the absolute value of $S$ by
\[
	\abs{S(x,t)} \lesssim \lim_{R\to\infty} \int_{-\pi/2}^{\pi/2}\exp\bigl(-(\,\abs{x}\Re l_{\alpha}(R\e^{\I\vphi})-t)R\cos\vphi\bigr)\dif \vphi.
\]
Let us write $\eps=\sqrt{\tau}\abs{x} - t >0$. For $R$ sufficiently large, $\abs{x}\Re l_{\alpha}(R\e^{\I\vphi})-t > \eps/2$, since $\Re l_{\alpha}(R\e^{\I\vphi}) \to \sqrt{\tau}$ by \eqref{eq: approx l_{alpha}}.  For such large $R$ the integrand is bounded by $\e^{-(\eps/2) R \cos\vphi}$, which converges pointwise to 0 and is bounded. From this it follows that the above integral converges to 0 when $R\to\infty$, by dominated convergence.

Suppose now that $\abs{x} < t/\sqrt{\tau}$. Then $\Re(-\abs{x}l_{\alpha}(s)+t) > 0$ for $\abs{s}$ sufficiently large. We will shift the contour to the left, resulting in a Hankel contour encircling the branch cut $(-\infty,0]$. For a small $\eps>0$, we set 
\begin{align*}
	\Gamma_{1}	&= [a-\I Y, -\I Y]; 						&	\Gamma_{5}	&= [\eps\e^{\I\pi}, Y\e^{\I\pi}]; \\
	\Gamma_{2}	&= \{Y\e^{\I\vphi}: -\pi/2\ge \vphi \ge -\pi\};		&	\Gamma_{6}	&=\{Y\e^{\I\vphi}: \pi\ge\vphi\ge\pi/2\};\\
	\Gamma_{3}	&= [Y\e^{-\I\pi}, \eps\e^{-\I\pi}]; 				&	\Gamma_{7}	&=[\I Y, a+ \I Y].\\
	\Gamma_{4}	&=\{\eps\e^{\I\vphi}: -\pi\le\vphi\le\pi\}; &&	
\end{align*}
By Cauchy's theorem, the contour integral in \eqref{eq: Laplace inverse} equals $\frac{1}{2\pi\I}\int_{\cup_{i}\Gamma_{i}}\tilde{S}(x,s)\e^{ts}\dif s$. On $\Gamma_{1}$ and $\Gamma_{7}$, $\tilde{S}(x,s)\e^{ts} \lesssim 1/Y$, so the integral over these pieces converges to zero as $Y\to\infty$. On $\Gamma_{2}$ and $\Gamma_{3}$, $\Re(-\abs{x}sl_{\alpha}(s) +ts) \sim (t-\sqrt{\tau}\abs{x})(\cos\vphi) Y$. As before, the integrals over these contours tend to zero by dominated convergence, since now $t-\sqrt{\tau}\abs{x} > 0$ and $\cos\vphi < 0$ (except at the boundary points $\vphi=\pm\pi/2$.) Since $\tilde{S}(x,s)\e^{ts} \sim 1/(2s)$ for $s\to0$, the integral over $\Gamma_{4}$ converges to $\I\pi$ as $\eps\to0$. Finally, combining the integrals over $\Gamma_{3}$ and $\Gamma_{5}$ and letting $Y\to\infty$, $\eps\to0$, we get the absolutely convergent integral in \eqref{eq: S in cone}. 
\end{proof}

\begin{remark}
(i) Proposition \ref{prop: support S} implies that the convolution in \eqref{eq: solution} is well-defined for arbitrary distributions $u_{0}$, $v_{0} \in \mathcal{D}'(\R)$,
and therefore the result given in \cite[Theorem 4.2]{KOZ2010} holds with $u_{0}$, $v_{0} \in \mathcal{D}'(\R)$. 

%(iii) \todo{Remark on boundary conditions}

(ii) In fact, the inverse Laplace integral \eqref{eq: Laplace inverse} converges absolutely for every $x$, $t$ with $x\neq0$. Writing $s=a+\I y$, the approximation \eqref{eq: approx l_{alpha}} shows that 
\begin{equation}\label{RealPartOfsl}
	\Re sl_{\alpha}(s) = \frac{\sqrt{\tau}}{2}\Bigl(\frac{1}{\tau}-1\Bigr)\sin(\alpha\pi/2)\abs{y}^{1-\alpha} + O(1+ \abs{y}^{1-2\alpha}), \quad \text{as } \abs{y}\to\infty,
\end{equation}
locally uniformly in $a$. Recall that $0<\alpha<1$ and $0<\tau<1$, so for non-zero $x$, the exponential in the integral decays like $\e^{-c\,\abs{x} \,\abs{y}^{1-\alpha}}$, where $c=\sqrt{\tau}/2(1/\tau-1)\sin(\alpha\pi/2)$ is a positive constant. This proves the absolute convergence. It is convenient to move the contour of integration to the line $\Re s=0$. For this we consider the following contours for a small parameter $\eps>0$:
\begin{align*}
	\Gamma_{1}	&= [\I Y, a+\I Y]; 						&\Gamma_{4}	&= [-\I Y, -\I\eps];\\
	\Gamma_{2}	&= [\I\eps, \I Y]; 						&\Gamma_{5}	&= [-\I Y, a-\I Y].\\
	\Gamma_{3}	&= \{\eps\e^{\I\vphi}: -\pi/2\le \vphi\le \pi/2\};	&&
\end{align*}
The estimate \eqref{RealPartOfsl} of $\Re sl_{\alpha}(s)$ immediately implies that the integrals of $\tilde{S}(x,s)\e^{ts}$ along the contours $\Gamma_{1}$ and $\Gamma_{5}$ tend to $0$ as $Y\to \infty$, whenever $x\neq0$. Since $\tilde{S}(x,s)\e^{ts} \sim 1/(2s)$ as $s\to0$, the integral along $\Gamma_{3}$ tends to $\I\pi/2$ as $\eps\to0$. Combining the integrals over $\Gamma_{2}$ and $\Gamma_{4}$ and letting $\eps\to0$, $Y\to\infty$ yields the following representation for $S(x,t)$ when $x\neq0$:
\begin{equation}
\label{eq: S for x!=0}
S(x,t) = \frac{1}{4} +  \frac{1}{4\pi\I}\pv\int_{-\infty}^{\infty}l_{\alpha}(\I y)\exp\bigl(-\abs{x}\I yl_{\alpha}(\I y) + \I ty\bigr)\frac{\dif y}{y},
\end{equation}
where $\pv$ denotes the Cauchy principal value.

The integral in \eqref{eq: S in cone} does not converge for $\abs{x}\ge t/\sqrt{\tau}$, so the representation \eqref{eq: S for x!=0} will be particularly useful for studying the behavior of $S$ near the boundary of the cone. On the other hand, the integral in \eqref{eq: S for x!=0} does not converge absolutely for $x=0$, so \eqref{eq: S in cone} will be useful for studying $S$ at small values of $x$.
\end{remark}

Combining both integral representations \eqref{eq: S in cone}, \eqref{eq: S for x!=0} of the fundamental solution $S$ allows us to give a complete description of its regularity. 

%%%%%%%%%%%%%%%%%%%%%%%%%%%%%%%%%%%%%%%%%%%%%%%%%%%%%%%%%%%%%%%%%%%%%%%%%%%%%%
%% MICRO-LOCAL ANALYSIS %%%%%%%%%%%%%%%%%%%%%%%%%%%%%%%%%%%%%%%%%%%%%%%%%%%%%%%%%%%%%%
%%%%%%%%%%%%%%%%%%%%%%%%%%%%%%%%%%%%%%%%%%%%%%%%%%%%%%%%%%%%%%%%%%%%%%%%%%%%%%
\section{Micro-local analysis of $S$}\label{MLA}
In this section, we provide the full micro-local analysis of the fundamental solution $S$ with respect to $C^{\infty}$ and $G^{\sigma}$, $\sigma\ge1$, extending previous results in the literature. In \cite{APS}, some regularity of $S$ was shown, namely that the map $t\mapsto S(x,t)$ is smooth for fixed $x\neq0$.
The micro-local analysis of $S$ was initiated in \cite{HOZ16}. It reaches a form of non-characteristic regularity for solutions to the Cauchy problem \eqref {eq: Cauchy problem Zener} (see \cite[Theorem 18.1.8]{Hormander}). Namely, \cite[Theorem 3.2]{HOZ16} states
\[
	\WF(u^+) \subseteq \Bigl\{(x,t;\xi,\eta): x\in\R, t > 0, \xi \neq 0, \eta = 0\, \text{ or } \,\eta^2 =\frac{1}{\tau} \xi^2 \Bigr\},
\]
where $u^+$ denotes the restriction of the solution to \eqref {eq: Cauchy problem Zener} to the forward time $t>0$. 
%We improve the above result by determining the wave front set of the fundamental solution $S$. 
This result suggests that apart from the frequencies $(\xi,0)$, also the frequencies orthogonal to the boundary of the light cone could be singular frequencies. We will show that this is not the case: $S$ is smooth on this boundary. This is in contrast with the classical wave equation, whose fundamental solution is singular at the boundary of the forward light cone, and for which the singular frequencies are those orthogonal to this boundary. 
%
%One could expect equality in the above formula, but 
%in contrast to the classical wave equation, the fundamental solution for the fractional wave equation is smooth on the boundary of the forward light cone. However, singularities on the cone become visible if one uses Gevrey notion of regularity instead of standard $C^\infty$. \todo{I am not perfectly happy with this but it is ok for to be sent to Andreas and Jasson}

\subsection{$C^\infty$-regularity} The following theorem provides the evaluation of the $C^{\infty}$-wave front set of $S$. In particular, it will imply that $S$ is smooth off the half line $x=0$, $t\ge0$.%The first theorem in this section shows that in above formula one can not put equality and that the points on the cone in the set on the right hand side are in fact "artefact". 
 \begin{theorem}
\label{thm: regularity S}
The fundamental solution $S$ is an $L^1_{\mathrm{loc}}$-function which is continuous on $\R^{2}\setminus \{(0,0)\}$. Its partial derivative with respect to $x$, $\partial_{x}S$, is discontinuous on the half-line $x=0$, $t>0$. Everywhere else, $S$ is of class $C^{\infty}$. In particular, for the wave front set we have
\[
	\WF(S) 	= \{(0,0;\xi,\eta): (\xi, \eta)\neq(0,0)\} \cup \{(0,t; \xi,0): t>0, \xi\neq0\}.
\]
\end{theorem}
\begin{proof}
The representations \eqref{eq: S in cone} and \eqref{eq: S for x!=0} imply the continuity of $S$ in the open sets $\abs{x} < t/\sqrt{\tau}$ and $x\neq0$ respectively, showing that $S$  coincides with a continuous function on $\R^{2}\setminus\{(0,0)\}$. It is not possible that $S$ contains linear combinations of $\delta^{(n)}(x)\delta^{(m)}(t)$, since these would show up in the Laplace transform $\tilde{S}$ as linear combinations of $s^{m}\delta^{(n)}(x)$, and they are not present in \eqref{eq: tilde{S}}. Hence, $S$ is an $L^{1}_{\mathrm{loc}}$-function, continuous on $\R^{2}\setminus\{(0,0)\}$.

Differentiating formula \eqref{eq: S in cone} with respect to $x$, and taking the limit for $x\to 0$ from the right and from the left, we see that
\begin{align*}
\dpd{S}{x}(0^{\pm}, t) 	&= \pm\frac{1}{4\pi\I}\int_{0}^{\infty}\biggl(\frac{1+\tau q^{\alpha}\e^{\I\alpha\pi}}{1+q^{\alpha}\e^{\I\alpha\pi}}
						-\frac{1+\tau q^{\alpha}\e^{-\I\alpha\pi}}{1+q^{\alpha}\e^{-\I\alpha\pi}}\biggr)\e^{-tq}\dif q \\
					&= \mp\frac{1}{2\pi}\int_{0}^{\infty}(1-\tau)\sin(\alpha\pi) \frac{q^{\alpha}}{1+2\cos(\alpha\pi) q^{\alpha} + q^{2\alpha}}\e^{-tq}\dif q.	
\end{align*}
The integrand of the last integral is positive when $q>0$, so the integral is non-zero. This shows that $\partial_{x}S(x,t)$ is not continuous at $x=0$. 

Using the representation \eqref{eq: S for x!=0}, we see that $S$ is smooth at points $(x,t)$ with $x\neq0$. Indeed, differentiating under the integral yields an additional factor which is of polynomial growth. Since the exponential in the integrand decays like $\lesssim \e^{-c\,\abs{x} \,\abs[0]{y}^{1-\alpha}}$, with $c=\sqrt{\tau}/2(1/\tau-1)\sin(\alpha\pi/2)$ (see \eqref{RealPartOfsl}), the integral remains convergent. Note that it is crucial here that $0<\alpha<1$ and $0<\tau<1$.

To compute the wave front set, we use \eqref{eq: S in cone}. Differentiating under the integral, we see that $\partial_{t}^{m}S(x,t)$ is bounded on compact subsets of the cone $\abs{x} < t/\sqrt{\tau}$, for each $m\in\N$. This implies that at a point $(0,t)$, $t>0$, the ``singular frequencies'' can only be along the positive and negative $(\xi,0)$-direction. Since $S$ is real-valued, its wave front set is symmetric about the origin in the frequency variables. Hence, both directions are present in the wave front set. Finally, since $PS(x,t) = \delta(x)\delta(t)$, and differential and convolution operators do not enlarge the wave front set, also $(0,0;\xi, \eta)\in \WF(S)$ for every $(\xi,\eta)\neq(0,0)$.
\end{proof}

From \eqref{eq: S in cone}, one readily sees that $S$ converges to $1/2$ for $t\to\infty$. Indeed, locally uniformly in $x$, the integral converges to $0$ as $t\to\infty$ by dominated convergence: for $q\le 1$, the integrand in \eqref{eq: S in cone} is dominated by $\e^{\,\abs{x}q}O(q^{\alpha-1} + \abs{x}q^{\alpha})$ (which can be seen by Taylor approximation using \eqref{eq: approx l_{alpha} 0}), while for $q\ge Q$, it is dominated by $\e^{(\,\abs{x}\sup l_{\alpha} - T)q}$ if $t\ge T$. 

We can also describe the behavior of $S$ for $(x,t)\to (0,0)$.
\begin{proposition}
Suppose $\lambda \in [0,1/\sqrt{\tau}]$. Then
\[
	\lim_{t\to0^{+}}S(\pm\lambda t, t) = \frac{\sqrt{\tau}}{2}(1-\lambda\sqrt{\tau}).
\]
%We can determine the nature of the singularity of $t\mapsto S(0,t)$ at $t=0$. We have that 
%\[
%	S(0,t) = \frac{\sqrt{\tau}}{2}H(t) + g(t),
%\]
%where $g$ is some continuous function supported on $\R_{\ge0}$, with $g(0)=0$.
\end{proposition}
\begin{proof}
The statement for $\lambda=1/\sqrt{\tau}$ is trivial, since $S(\pm t/\sqrt{\tau}, t)=0$ for $t>0$.

Suppose first that $\lambda=0$. We compute the Laplace transform of $S(0,t) - (\sqrt{\tau}/2)H(t)$. We get
\[
	\ML\Bigl(S(0,t) - \frac{\sqrt{\tau}}{2}H(t)\Bigr)(s) = \tilde{S}(0,s) - \frac{\sqrt{\tau}}{2s} = \frac{l_{\alpha}(s)-\sqrt{\tau}}{2s}.
\]
By \eqref{eq: approx l_{alpha}}, $l_{\alpha}(s) = \sqrt{\tau}+O(\,\abs{s}^{-\alpha})$ as $\abs{s} \to \infty$. Therefore, the above Laplace transform decays as $\lesssim \abs{s}^{-1-\alpha}$, so it is integrable on every vertical line $\Re s = a$, with $a>0$. Hence, $S(0,t)-(\sqrt{\tau}/2)H(t)$ is a continuous function, which vanishes for $t\le0$.

For general $\lambda\in(0,1/\sqrt{\tau})$, we apply the same strategy, namely determining the asymptotic behavior of the Laplace transform of $t\mapsto S(\lambda t, t)$. Let $s=2+\I y$. Using the Laplace inversion for $S(x,t)$ \eqref{eq: Laplace inverse} with $a=1$, we get 
\[
	\ML\{S(\lambda t, t)\}(s) = \frac{1}{4\pi\I}\int_{0}^{\infty}\e^{-st}\int_{1-\I\infty}^{1+\I\infty}\frac{l_{\alpha}(z)}{z}\exp\bigl(tz(1- \lambda l_{\alpha}(z))\bigr)\dif z \dif t.
\]
We want to swap the order of integration here, using the Fubini-Tonelli theorem. This is allowed, since
\begin{multline*}
	\int_{1-\I\infty}^{1+\I\infty}\int_{0}^{\infty}\abs{\frac{l_{\alpha}(z)}{z}\exp\bigl(t(z-\lambda z l_{\alpha}(z) - s)\bigr)}\dif t \abs{\dif z} \\
	= \int_{1-\I\infty}^{1+\I\infty}\abs{\frac{l_{\alpha}(z)}{z}}\frac{1}{-1+\lambda \Re (zl_{\alpha}(z)) + 2} \abs{\dif z} < \infty
\end{multline*}
Here we used that $1-\lambda\Re(zl_{\alpha}(z)) - 2 \le -1$ and that $\lambda\Re(zl_{\alpha}(z)) \sim -\lambda c \abs{\Im z}^{1-\alpha}$ for some $c>0$, by \eqref{RealPartOfsl}, so that the above integral converges absolutely. Swapping the order of integration and integrating with respect to $t$, we get
\[
	\ML\{S(\lambda t, t)\}(s) = -\frac{1}{4\pi\I}\int_{1-\I\infty}^{1+\I\infty}\frac{l_{\alpha}(z)}{z\bigl(z(1-\lambda l_{\alpha}(z)) - s\bigr)}\dif z.
\]
We will evaluate this integral via Cauchy's theorem. Let $y=\Im s$ be large but fixed. The integrand above decays like $\lesssim 1/\abs{z}^{2}$, and has a unique pole in the right half plane $\Re z > 1$, which we denote by $z(s)$ (this follows for example by applying the argument principle on the region enclosed by the line $[1-\I R, 1+\I R]$ and the right semicircle with center $1$ and radius $\sqrt{1+R^{2}}$, for sufficiently large $R$). We have
\[
	\bigl(1-\lambda l_{\alpha}(z(s))\bigr) z(s) - s = 0, \quad z(s) = \frac{s}{1-\lambda\sqrt{\tau}}\bigl(1+O(\,\abs{s}^{-\alpha})\bigr),
\] 
where we used \eqref{eq: approx l_{alpha}}.
Applying Cauchy's theorem and \eqref{eq: approx l_{alpha}} again, we get
\[
	\ML\{S(\lambda t, t)\}(s) = \frac{l_{\alpha}(z(s))}{2z(s)} = \frac{\sqrt{\tau}(1-\lambda\sqrt{\tau})}{2s}\bigl(1+O(\,\abs{s}^{-\alpha})\bigr), \quad \text{ for } \abs{y}=\abs{\Im s} \to \infty.
\]
It follows that $S(\lambda t, t) - (\sqrt{\tau}/2)(1-\lambda\sqrt{\tau})H(t)$ is a continuous function, vanishing for $t\le0$, since its Laplace transform, decaying like $\lesssim \abs{s}^{-1-\alpha}$, is absolutely integrable on the line $\Re s = 2$.
\end{proof}

\subsection{Gevrey regularity} Using a finer notion of smoothness, namely by means of the Gevrey classes $G^{\sigma}$, more singularities become visible. In the following theorem, we describe the $G^{\sigma}$-regularity of $S$ for every $\sigma \in [1, \infty)$. We see that for $\sigma$ sufficiently close to $1$, namely $1\le \sigma < 1/(1-\alpha)$, the boundary of the forward light cone becomes singular.
\begin{theorem}
\label{th: Gevrey regularity}
On $\R^{2}\setminus \bigl(\{0\}\times[0,\infty)\bigr)$, $S$ belongs to the Gevrey class $G^{\frac{1}{1-\alpha}}$. Furthermore, at points $(x,t)$ with $\abs{x}\neq t/\sqrt{\tau}$ and $x\neq0$ it is real analytic. For the wave front set with respect to $G^{\sigma}$, we have the following equalities:
\begin{align*}
	\WF_{G^{\sigma}}(S) 	= {}	&\{(0,0;\xi,\eta): (\xi, \eta)\neq(0,0)\} \cup \{(0,t; \xi,0): t>0, \xi\neq0\}		&&\text{if } \sigma\ge\frac{1}{1-\alpha}; \\
	\WF_{G^{\sigma}}(S)	= {}	&\{(0,0;\xi,\eta): (\xi, \eta)\neq(0,0)\} \cup\{(0,t; \xi,0): t>0, \xi\neq0\} &&\\
				{}\cup{}	&\{(x, t; \xi, \eta): t>0, (\xi,\eta)\neq(0,0),\ \abs{x}=t/\sqrt{\tau}, (x,t)\cdot(\xi,\eta)=0\}	&&\text{if }1 \le \sigma < \frac{1}{1-\alpha}.
\end{align*}
%\[ 
%	\WF_{G^{\sigma}}(S) = \begin{dcases}
%						\{(0,0,\xi,\eta): (\xi, \eta)\neq(0,0)\} \cup \{(0,t; \xi,0): t>0, \xi\neq0\}		&\text{if } \sigma\ge\frac{1}{1-\alpha}; \\
%						\begin{split}
%							&\{(0,0,\xi,\eta): (\xi, \eta)\neq(0,0)\} \cup \{(0,t; \xi,0): t>0, \xi\neq0\} \\
%							&{} \cup\{(x, t; \xi, \eta): t>0, (\xi,\eta)\neq(0,0),\ \abs{x}=t/\sqrt{\tau}, (x,t)\cdot(\xi,\eta)=0\}
%						\end{split}
%																				&\text{if }1 \le \sigma < \frac{1}{1-\alpha}.		
%	\end{dcases}
%\]
\end{theorem}

\begin{proof}
The representation \eqref{eq: S in cone} readily implies that $S$ is real analytic at points $(x,t)$ with $x\neq0$, $\abs{x} < t/\sqrt{\tau}$. Indeed, the integral and its derivates with respect to $x$ and $t$ still converge when one replaces $(x,t)$ by $(x+z_{1}, t+z_{2})$, $z_{1}$, $z_{2} \in \C$ with $\abs{z_{1}}$ and $\abs{z_{2}}$ sufficiently small. 
%In the same way we see that the function $t \mapsto S(0,t)$ is a real analytic function of $t$ for $t>0$.

Let us now see that $S$ is in the Gevrey class $G^{\frac{1}{1-\alpha}}$ for $x\neq0$. Suppose that $x>0$. Differentiating under the integral sign in \eqref{eq: S for x!=0}, we see that 
\[
	\pd[n]{}{x}\pd[m]{}{t}S(x,t) = \frac{(-1)^{n}}{4\pi\I}\int_{-\infty}^{\infty}\bigl(l_{\alpha}(\I y)\bigr)^{n}(\I y)^{n+m}\exp\bigl(-\I xy l_{\alpha}(\I y) + t\I y\bigr) \frac{\dif y}{y}.
\]
Using \eqref{eq: bound l_{alpha}}, this is in absolute value bounded by 
\[
	D_{1}^{n+1}\biggl(\int_{0}^{1}y^{n+m-1}\e^{-c_{1}xy^{1+\alpha}}\dif y + \int_{1}^{\infty}y^{n+m-1}\e^{-c_{2}xy^{1-\alpha}}\dif y\biggr) 
	\le D_{2}^{n+m+1}\Bigl\{1+x^{-\frac{n+m}{1-\alpha}}\Gamma\Bigl(\frac{n+m}{1-\alpha}\Bigr)\Bigr\},
	%\exp\biggl\{-x\biggl(\sin(\alpha\pi/2) \frac{\sqrt{\tau}}{2}\biggl(\frac{1}{\tau}-1\biggr)y^{1-\alpha} + O(y^{1-2\alpha})\biggr)\biggr\}\dif y.
\] 
for some positive constants $D_{1}, D_{2}$, by changing variables $y' = c_{2}xy^{1-\alpha}$ in the second integral. For $x$ in a closed subset $F$ of $\R\setminus\{0\}$, this is bounded by $D_{F}^{n+m+1}(n!m!)^{\frac{1}{1-\alpha}}$, where $D_{F}$ is a positive constant depending on $F$.
%A routine calculation now shows that this is bounded by $C_{2}^{n+m+1}\bigl\{1+x^{-\frac{n+m}{1-\alpha}}\Gamma\bigl(\frac{n+m}{1-\alpha}\bigr)\bigr\}$
%((n+m)!)^{\frac{1}{1-\alpha}}\le C_{2}^{n+m+1}2^{\frac{n+m}{1-\alpha}}(n!m!)^{\frac{1}{1-\alpha}}$, where $C_{2}$ is some positive constant depending on $\delta$.

Next, we will show that for $\sigma<1/(1-\alpha)$, the boundary of the forward light cone is contained in $\singsupp_{G^{\sigma}} S$. For this it is convenient to perform the change of variables $u=\sqrt{\tau} x + t$, $v=\sqrt{\tau} x - t$, so that points with $(x,t)$ satisfying $x>0$ and $x = t/\sqrt{\tau}$ have new coordinates $(u, 0)$ with $u>0$. (To treat the boundary points with $x<0$, one considers an analogous change of variables, or one uses the symmetry $S(-x,t) = S(x,t)$.) We set 
\[
	S^{\natural}(u,v) = S(x, t) = S\Bigl(\frac{u+v}{2\sqrt{\tau}}, \frac{u-v}{2}\Bigr).
\]

Let $u>0$ and $\sigma<1/(1-\alpha)$. We claim that $(u,0) \in \singsupp_{G^{\sigma}} S^{\natural}$. This is equivalent to the statement that for every neighborhood $U$ of $(u,0)$ and for every $C>0$ there exists some $\beta\in \N^{2}$, and some point $(a,b)\in U$ for which
\[
	\abs{\partial^{\beta}S^{\natural}(a,b)} > C^{1+\,\abs{\beta}}(\beta!)^{\sigma}.
\]
We will actually prove something stronger, namely that there exists a constant $C>0$ and a sequence $(v_{m})_{m}$, $v_{m}>0$, $v_{m}\to 0$ (both $C$ and $v_{m}$ depending on $u$) so that 
\begin{equation}
\label{eq: Gevrey lower bound}
	\abs[3]{\dpd[m]{S^{\natural}}{v}(u,-v_{m})} \gtrsim C^{m}(m^{m})^{\frac{1}{1-\alpha}},
\end{equation}
provided that $m$ is sufficiently large. Showing this estimate requires an intricate technical analysis of the inverse Laplace integral. 
%We first give an overview of the method for showing \eqref{eq: Gevrey lower bound}, before delving into the more technical details. 
First we write $\partial^{m} S^{\natural}/\partial v^{m}$ as a contour integral:
\[
	\pd[m]{S^{\natural}}{v}(u,v) = \frac{(-1)^{m}}{2\pi}\Im \int_{0}^{+\I\infty}l_{\alpha}(s)\Bigl(\frac{s}{2}\Bigr)^{m}\Bigl(\frac{l_{\alpha}(s)}{\sqrt{\tau}}+1\Bigr)^{m}
							\exp\biggl(-\frac{us}{2}\Bigl(\frac{l_{\alpha}(s)}{\sqrt{\tau}}-1\Bigr) - \frac{vs}{2}\Bigl(\frac{l_{\alpha}(s)}{\sqrt{\tau}}+1\Bigr)\biggr)\frac{\dif s}{s}.
\]
Here, we integrate along the half-line $[0, +\I\infty)$, and we used the symmetry $\int_{-\I\infty}^{0} = -\overline{\int_{0}^{+\I\infty}}$ to write the inverse Laplace integral as the imaginary part of the integral over the part of the contour with $\Im s \ge 0$.

The first step is to perform a change of variables in the above integral, which will bring out the $(m^{m})^{\frac{1}{1-\alpha}}$-behavior, and to substitute a well-chosen value for $v$, which will in some sense simplify the remaining integral. Namely, we will set
\[
	s = \mu m^{\frac{1}{1-\alpha}} w, \quad \text{ with } \mu = \biggl(\frac{u}{4}\biggl(\frac{1}{\tau}-1\biggr)\frac{1}{1-\alpha}\biggr)^{-\frac{1}{1-\alpha}}, 
	\quad \text{ and } v = -v_{m} = -\frac{\kappa}{\mu}m^{-\frac{\alpha}{1-\alpha}}.
\]
Here, $\kappa$ is some large number depending on $m$, whose value will be chosen later. For the moment, we only specify a fixed range for $\kappa$, say\footnote{The value $1000$ occurring here is somewhat arbitrary, we just require some large fixed number.} $1000/\sin(\alpha\pi) \le \kappa^{1-\alpha} \le 2000/\sin(\alpha\pi)$.

With the above substitution we get
\begin{align*}
	\dpd[m]{S^{\natural}}{v}(u,-v_{m}) = \frac{(-1)^{m}}{2\pi}\Bigl(\frac{\mu}{2}\Bigr)^{m}(m^{m})^{\frac{1}{1-\alpha}}\Im \int_{0}^{+\I\infty}l_{\alpha}(\mu m^{\frac{1}{1-\alpha}}w)w^{m}\biggl(\frac{l_{\alpha}(\mu m^{\frac{1}{1-\alpha}}w)}{\sqrt{\tau}} + 1\biggr)^{m}\\
	\exp\biggl\{-\frac{u\mu m^{\frac{1}{1-\alpha}}w}{2}\biggl(\frac{l_{\alpha}(\mu m^{\frac{1}{1-\alpha}}w)}{\sqrt{\tau}} - 1\biggr) + \frac{m\kappa w}{2}\biggl(\frac{l_{\alpha}(\mu m^{\frac{1}{1-\alpha}}w)}{\sqrt{\tau}} + 1\biggr)\biggr\}\frac{\dif w}{w}.
\end{align*}
Let us now consider the remainders
\begin{equation}
\label{eq: remainders}
	E_{1}(s) =  \frac{l_{\alpha}(s)}{\sqrt{\tau}} - 1, \quad E_{2}(s) = \frac{l_{\alpha}(s)}{\sqrt{\tau}} - 1 - \frac{1}{2}\biggl(\frac{1}{\tau}-1\biggr)s^{-\alpha},
\end{equation}
then by \eqref{eq: approx l_{alpha}} we have 
\begin{equation}
\label{eq: bounds E1 and E2}
E_{1}(s) \lesssim \abs{s}^{-\alpha}, \quad  \text{and} \quad E_{2}(s) \lesssim \abs{s}^{-2\alpha}, \quad  \text{as } \abs{s} \to \infty.
\end{equation}
The expression for $\partial_{v}^{m}S^{\natural}(u, -v_{m})$ can be rewritten as
\begin{align*}
	\frac{(-1)^{m}}{2\pi}\Bigl(\frac{\mu}{2}\Bigr)^{m}(m^{m})^{\frac{1}{1-\alpha}}\Im \int_{0}^{+\I\infty}l_{\alpha}(\mu m^{\frac{1}{1-\alpha}}w)
	\exp\biggl\{m\biggl(\kappa w - (1-\alpha)w^{1-\alpha} + \log w+\log 2 \\
	+ \biggl[ \frac{\kappa w}{2}E_{1}(\mu m^{\frac{1}{1-\alpha}} w) - \frac{u\mu m^{\frac{\alpha}{1-\alpha}}w}{2}E_{2}(\mu m^{\frac{1}{1-\alpha}} w) + \log\Bigl(1+\frac{E_{1}(\mu m^{\frac{1}{1-\alpha}} w)}{2}\Bigr)\biggr] \biggr)\biggr\}\frac{\dif w}{w}.
\end{align*}
Denoting the terms in between the square brackets $[ \dotso ]$ by $g_{m}(w)$, and setting 
\[
	f(w) \coloneqq \kappa w - \frac{1}{1-\alpha}w^{1-\alpha} + \log w+ \log 2,
\] 
we have
\[
	\pd[m]{S^{\natural}}{v}(u,-v_{m}) = \frac{(-1)^{m}}{2\pi}\Bigl(\frac{\mu}{2}\Bigr)^{m}(m^{m})^{\frac{1}{1-\alpha}}\Im \int_{0}^{+\I\infty}l_{\alpha}(\mu m^{\frac{1}{1-\alpha}}w) \exp\bigl(m(f(w) + g_{m}(w))\bigr)\frac{\dif w}{w}.
\]

\begin{lemma}
\label{lem: technical lemma}
For $m$ sufficiently large, one can choose $\kappa= \kappa_{m}$ in the fixed range $1000/\sin(\alpha\pi) \le \kappa^{1-\alpha} \le 2000/\sin(\alpha\pi)$ such that 
\begin{equation*}
\Im \int_{0}^{+\I\infty}l_{\alpha}(\mu m^{\frac{1}{1-\alpha}}w)\exp\bigl(m(f(w) + g_{m}(w))\bigr)\frac{\dif w}{w} \gtrsim \frac{c^{m}}{\sqrt{m}}. 
\end{equation*}
Here, $c$ is a positive constant independent of $m$.
\end{lemma}

We will not give the proof here, since it is rather lengthy and technical. Instead we provide a proof in Appendix \ref{sec: technical calculation}. The main idea is the following. In view of the estimates \eqref{eq: bounds E1 and E2}, $g$ is small for large $m$. Also $l_{\alpha}(\mu m^{\frac{1}{1-\alpha}} w) \to \sqrt{\tau}$ as $m\to\infty$. In some sense, the analysis reduces to the analysis of the simpler integral $\int (\sqrt{\tau}/w)\e^{mf(w)}\dif w$, which can be estimated with the saddle point method. Indeed, the constant $c$ is related to the value of $\e^{\Re f}$ at the saddle point $w_{0}$ of $f$. The purpose of the free parameter $\kappa$ is to control the imaginary part of the integral. For each $m$, we will choose a $\kappa_{m}$ so that the argument of $\int (\sqrt{\tau}/w)\e^{mf(w)}\dif w$ is close to $\pi/2$.

Assuming Lemma \ref{lem: technical lemma}, we get that
\[
	\abs[3]{\dpd[m]{S^{\natural}}{v}(u,-v_{m})} \gtrsim \frac{1}{\sqrt{m}}\Bigl(\frac{\mu c}{2}\Bigr)^{m}(m^{m})^{\frac{1}{1-\alpha}},
\]	 
from which \eqref{eq: Gevrey lower bound} follows for any $C < \mu c/2$. This shows that $(u,0)\in \singsupp_{G^{\sigma}} S^{\natural}$, as soon as $\sigma < 1/(1-\alpha)$.

To determine the wave front set with respect to $G^{\sigma}$, we show that at the point $(u_{0},0)$, $u_{0}>0$, $S^{\natural}$ is ``real analytic in the $(u,0)$-direction''. Indeed, taking partial derivatives with respect to $u$ we get from \eqref{eq: S for x!=0} with $v \ge -u/2$
\[
	\pd[n]{S^{\natural}}{u}(u,v) \lesssim \int_{0}^{\infty} \Bigl(\frac{y}{2}\Bigr)^{n}\abs{\frac{l_{\alpha}(\I y)}{\sqrt{\tau}} - 1}^{n}\exp\Bigl(\frac{u}{4\sqrt{\tau}}y \Im l_{\alpha}(\I y) \Bigr)\frac{\dif y}{y}.
\]
Now by \eqref{eq: approx l_{alpha}}, $\abs{l_{\alpha}(\I y)/\sqrt{\tau} - 1} \lesssim y^{-\alpha}$ for large $y$. Using this and \eqref{eq: bound l_{alpha}}, we get 
\begin{align*}
	\dpd[n]{S^{\natural}}{u}(u,v) 
		&\lesssim D_{1}^{n}\biggl(\int_{0}^{1}y^{n-1}\exp\Bigl(-\frac{uc_{1}}{4\sqrt{\tau}}y^{1+\alpha}\Bigr)\dif y +
		\int_{1}^{\infty}(y^{1-\alpha})^{n}\exp\Bigl(-\frac{uc_{2}}{4\sqrt{\tau}}y^{1-\alpha}\Bigr)\frac{\dif y}{y}\biggr) \\
		&\le D_{1}^{n}(1+ D_{2}^{n}n!).
	%\quad C_{2} = C_{1} \cdot \frac{4\sqrt{\tau}}{Au}.
\end{align*}
Here, $D_{2} = D_{2}(u) = 4\sqrt{\tau}/(uc_{2}(1-\alpha))$ can be bounded uniformly on some neighborhood of $(u_{0},0)$. This implies that at the point $(u_{0},0)$, the ``$G^{1}$-singular frequencies\footnote{i.e.\ the points $(\xi,\eta)$ with $(u_{0},0;\xi,\eta)\in\WF_{G^{1}}(S^{\natural})$}\,'' can only occur along the positive and negative $(0,\eta)$-direction. Since $WF_{G^{\sigma}}(S^{\natural})$ is symmetric about the origin in the frequency variables, and since $(u_{0},0) \in \singsupp_{G^{\sigma}} S^{\natural}$ for any $\sigma<1/(1-\alpha)$, we get for such $\sigma$ that $(u_{0},0;\xi, \eta) \in WF_{G^{\sigma}}(S^{\natural}) \iff \xi=0$ and $\eta\neq0$. 

A similar argument, now using representation \eqref{eq: S in cone}, shows that at points $(0,t_{0})$ with $t_{0}>0$, $S$ is real analytic in the $(0,t)$-direction, so that the wave front set with respect to $G^{\sigma}$, $\sigma\ge1$, can only contain directions orthogonal to the line $x=0$ at points $(0,t_{0})$. This completes the proof of the theorem.
\end{proof}

%%%%%%%%%%%%%%%%%%%%%%%%%%%%%%%%%%%%%%%%%%%%%%%%%%%%%%%%%%%%%%%%%%%%%%%%%%%%%%
%% QUALITATIVE ANALYSIS %%%%%%%%%%%%%%%%%%%%%%%%%%%%%%%%%%%%%%%%%%%%%%%%%%%%%%%%%%%%%%%
%%%%%%%%%%%%%%%%%%%%%%%%%%%%%%%%%%%%%%%%%%%%%%%%%%%%%%%%%%%%%%%%%%%%%%%%%%%%%%
\section{Qualitative analysis}\label{QA}
In this section, we will discuss some qualitative aspects of the Fractional Zener wave equation and some of its solutions. First, we consider so-called pseudo-monochromatic waves as a means to study dispersion and dissipation. Next, we will analyse the solution of the Cauchy problem \eqref{eq: Cauchy problem Zener} with initial condition a delta concentrated at the origin. Studying this solution will allow us to define a meaningful notion of ``wave speed'' for this equation.

%%%%%%%%%%%%%%%%%%%%%%%%%%%%%%%%%%%%%%%%%%%%%%%%%%%%%%%%%%%%%%%%%%%%%%%%%%%%%%
\subsection{Dispersion and dissipation}
When studying dispersion, one investigates the relation between the phase velocity $V(\omega)$ and the frequency $\omega$ of a wave solution $A\e^{\I\omega(t-x/V(\omega))}$. In the absence of such purely monochromatic wave solutions of the homogeneous Cauchy problem, we will investigate the response when we submit the system to a forced oscillation with frequency $\omega$: let $u$ be the solution of \eqref{eq: Cauchy problem Zener} with initial conditions $u_{0}=v_{0}=0$ and force term $f(x,t) = \delta(x)H(t)\cos(\omega t)$ for some $\omega>0$. Let us first mention, for the sake of comparison, the solution $u_{\mathrm{cl}}$ to the \emph{classical} wave equation with these Cauchy data:
\begin{gather*}
	\biggl(\dpd[2]{}{t} - \frac{1}{\tau}\dpd[2]{}{x}\biggr)u_{\mathrm{cl}}(x,t) = f(x,t), \quad f(x,t)= \delta(x)H(t)\cos(\omega t), \quad u_{0}(x)=v_{0}(x)=0 \\
	\implies u_{\mathrm{cl}}(x,t) = H(t/\sqrt{\tau} - \abs{x})\frac{\sqrt{\tau}}{2\omega}\sin\bigl(\omega t - \sqrt{\tau}\omega \abs{x}\bigr).
\end{gather*}
This solution represents two waves traveling in opposite directions. They have wave number $k$ related to the frequency $\omega$ via the simple dispersion relation $k(\omega) = \sqrt{\tau}\omega$, and phase speed $V(\omega) = 1/\sqrt{\tau}$. 

Let us now analyze the solution in the fractional Zener case.
In view of Theorem \ref{thm: regularity S}, the solution $u=S\ast f$ is smooth for $x\neq0$. It has Laplace transform
\[
	\tilde{u}(x,s) = \frac{l_{\alpha}(s)}{2}\e^{-\,\abs{x}sl_{\alpha}(s)}\frac{1}{s^{2}+\omega^{2}}.
\]
From Proposition \ref{prop: support S} it follows that $u(x,t)=0$ if $\abs{x}>t/\sqrt{\tau}$. If $\abs{x} < t/\sqrt{\tau}$, we transform the contour to the contour which encircles the negative real axis, like was done to deduce \eqref{eq: S in cone}. However in this case, we get two contributions from the poles at $s=\pm\I\omega$, and no contribution from $s=0$. We get
\[
	u(x,t) = H(t/\sqrt{\tau} - \abs{x}) \bigl(u_{\mathrm{ss}}(x,t) + u_{\mathrm{ts}}(x,t)\bigr),
\]
where, using the notation $l_{\alpha}(\I\omega) = a(\omega)-\I b(\omega) = \rho(\omega)\e^{-\I\phi(\omega)}$ with $\sgn b(\omega) = \sgn \phi(\omega) = \sgn\omega$, 
\begin{align}
	u_{\mathrm{ss}}(x,t)	&= \frac{l_{\alpha}(\I\omega)}{4\I\omega}\e^{-\,\abs{x}\I\omega l_{\alpha}(\I\omega) + \I\omega t} - \frac{l_{\alpha}(-\I\omega)}{4\I\omega}\e^{\,\abs{x}\I\omega l_{\alpha}(-\I\omega) - \I\omega t} \nonumber\\
						&=\frac{\rho(\omega)}{2\omega}\e^{-b(\omega)\omega\,\abs{x}}\sin\bigl(\omega t - a(\omega)\omega \abs{x} - \phi(\omega)\bigr); \label{eq: steady state}
\end{align} 
\begin{equation}						
	u_{\mathrm{ts}}(x,t)		= \frac{1}{4\pi\I}\int_{0}^{\infty}\bigl(l_{\alpha}(q\e^{-\I\pi})\e^{\,\abs{x}l_{\alpha}(q\e^{-\I\pi})q} - l_{\alpha}(q\e^{\I\pi})\e^{\,\abs{x}l_{\alpha}(q\e^{\I\pi})q}\bigr)\frac{\e^{-tq}}{q^{2}+\omega^{2}}\dif q. \nonumber
\end{equation} 
We call $u_{\mathrm{ss}}$ the steady state, and $u_{\mathrm{ts}}$ the transient state. Indeed, from the above formula it is clear that $u_{\mathrm{ts}}(x,t) \to 0$ as $t\to\infty$, locally uniformly in $x$.

%This can be compared with the classical wave equation with wave speed $1/\sqrt{\tau}$, where the forced oscillation $f(x,t)=\delta(x)H(t)\cos(\omega t)$ yields the response

Following Mainardi \cite[Section 4.3]{Mainardi}, we call the steady state \eqref{eq: steady state} a ``pseudo-monochromatic wave'' with complex wave number $k$ satisfying the dispersion relation $k(\omega)=\omega l_{\alpha}(\I\omega)$.
%and complex refraction index $n(\omega) = l_{\alpha}(\I\omega)/\sqrt{\tau}$. 
It has phase velocity \[V(\omega) = 1/a(\omega),\] and has an amplitude which is exponentially decreasing in space, with attenuation coefficient $d(\omega) = b(\omega)\omega$.  The exponential dampening in space indicates dissipation, quantified by the attenuation coefficient, which has the following asymptotics, following from Lemma \ref{lem: l_{alpha}}:
\begin{align}
	d(\omega)		&\sim \frac{1-\tau}{2}\sin(\alpha\pi/2)\omega^{1+\alpha}, 	 			&&\omega	\to 0,	\nonumber\\	
	d(\omega) 	&\sim \frac{\sqrt{\tau}}{2}(1/\tau - 1)\sin(\alpha\pi/2)\omega^{1-\alpha}, 	&&\omega	\to\infty.  \label{eq: asymp behavior delta(omega)}
\end{align}
Since $V(\omega)$ is non-constant, there is some dispersion; however $V(\omega)$ is nearly constant, in the sense that it increases monotonically from $1$ to $1/\sqrt{\tau}$ when $\omega$ increases from $0$ to $\infty$.
The fact that the phase velocity $V(\omega)$ is increasing in $\omega$, indicates that the dispersion is anomalous. One may define the group velocity as 
\[
	U(\omega) = \Bigl(\od{(\Re k)}{\omega}(\omega)\Bigr)^{-1} = \frac{1}{a(\omega) + \omega a'(\omega)}.
\]
Note that $U(\omega) \ge V(\omega)$, with equality only for $\omega=0$ and in the limit $\omega\to\infty$. In the presence of dissipation and anomalous dispersion, this notion of group velocity looses its physical interpretation as velocity of a wave packet. However, in the next subsection we will provide a natural definition for the wave packet speed.

%%%%%%%%%%%%%%%%%%%%%%%%%%%%%%%%%%%%%%%%%%%%%%%%%%%%%%%%%%%%%%%%%%%%%%%%%%%%%%
\subsection{Shape of the wave packet}
%In view of Proposition \ref{prop: support S}, one might be tempted to define the wave speed corresponding to the fractional Zener model as simply $1/\sqrt{\tau}$. For the classical wave equation 
%\[
%	\biggl(\dpd[2]{}{t} - \frac{1}{\tau}\dpd[2]{}{x}\biggr)u = 0,
%\]
%this answer is satisfactory, since in this case the solution of the homogeneous problem with initial condition $u_{0}(x)$ and $v_{0}(x)=0$ is given by 
%\[
%	u(x,t) = \frac{1}{2}u_{0}(x-t/\sqrt{\tau}) + \frac{1}{2}u_{0}(x+t/\sqrt{\tau}),
%\]
%which clearly represents two waves moving in opposite directions with speed $1/\sqrt{\tau}$. In the fractional Zener model, the situation is less straightforward. We will show that, although the wave front moves with speed $1/\sqrt{\tau}$, the ``bulk'' of the wave moves at a slower speed. We therefore find it more appropriate to refer to $1/\sqrt{\tau}$ as the wave front speed, and reserve the term ``wave speed'' for this slower speed, which turns out to be equal to 1.

We denote by $K(x,t)$ the solution of \eqref{eq: Cauchy problem Zener} with Cauchy data $u_{0}(x) = \delta(x)$, $v_{0}(x) =0$, $f(x,t)=0$; so $K(x,t)=\partial_{t}S(x,t)$. %In \cite{KOZ2010}, $K$ was plotted for some values of the parameters $\alpha$ and $\tau$, and for some 
In this section we will accurately describe the shape of the ``wave packet'' $K(x,t)$, when $t$ is sufficiently large. 
%We will see that, although the wave front speed of this wave packet is given by $1/\sqrt{\tau}$, the packet itself moves at the smaller speed $1$. %\todo{This solution was plotted in \cite{KOZ2010} indicated }

The solution with general initial condition $u_{0}(x)$ and $v_{0}(x) = 0$ is given by $u(x,t)=K(x,t)\ast_{x} u_{0}(x)$. The evolution of a general wave packet with initial shape $u_{0}(x)$ can then be described using the analysis of $K$. Let us first list some simple properties of $K$.
\begin{proposition}
\label{prop: basic properties K}
For any fixed $t>0$, the function $x \mapsto K(x,t)$ is an even continuous function of $x$, supported in $[-t/\sqrt{\tau}, t/\sqrt{\tau}]$, with integral $1$. %\todo{Positivity?} 
\end{proposition}
\begin{proof}
The continuity, evenness, and statement on the support all follow from properties of $S$. In order to compute the integral $\int_{-\infty}^{\infty}K(x,t)\dif x$, we consider the representation of $K$ as inverse Laplace transform:
\begin{equation}
\label{eq: K}
	K(x,t) =\frac{1}{4\pi\I}\int_{a-\I \infty}^{a+\I \infty}l_{\alpha}(s)\exp\bigl(-\abs{x}s l_{\alpha}(s) + st\bigr)\dif s, \quad x\neq0.
\end{equation}
We have
\begin{align*}
	\int_{-\infty}^{\infty}K(x,t)\dif x 	&= 2\lim_{\eps\to0^{+}}\int_{\eps}^{\infty}K(x,t)\dif x = \lim_{\eps\to0^{+}}\frac{1}{2\pi\I}\int_{a-\I\infty}^{a+\I\infty}\frac{1}{s}\exp\bigl(-\eps s l_{\alpha}(s) +st\bigr)\dif s = 1.
\end{align*}
Here we introduced the parameter $\eps$ to be able to swap the order of integration, which is allowed for $\eps>0$ by the Fubini-Tonelli theorem. The last equality follows for example by shifting the contour to a Hankel contour, as in \eqref{eq: S in cone}: if $\eps< t/\sqrt{\tau}$, one may write the above integral as
\[
	\frac{1}{2\pi\I}\int_{a-\I\infty}^{a+\I\infty}\frac{1}{s}\exp\bigl(-\eps s l_{\alpha}(s) +st\bigr)\dif s
	 = 1 + \frac{1}{2\pi\I}\int_{0}^{\infty}\bigl\{\exp\bigl(\eps ql_{\alpha}(q\e^{\I\pi})\bigr) - \exp\bigl(\eps ql_{\alpha}(q\e^{-\I\pi})\bigr)\bigr\}\frac{\e^{-tq}}{q}\dif q.
\]
The last integral converges to $0$ as $\eps\to0$ by dominated convergence. Indeed, the integrand converges pointwise to $0$, and for the dominating function, we can argue as follows. Suppose that $\eps < t/(2\sqrt{\tau})$, and let $Q$ be such that $q\ge Q \implies \Re l_{\alpha}(q\e^{\pm\I\pi}) \le 4\sqrt{\tau}/3$. For $0\le q \le Q$, we apply Taylor's theorem to see that the integrand is dominated by $O_{Q}(\e^{-tq})$. For $q\ge Q$, the integrand is dominated by $\e^{-tq/3}$.
\end{proof}

 \begin{figure}[htb]
 \centerline{\includegraphics[width=7cm]{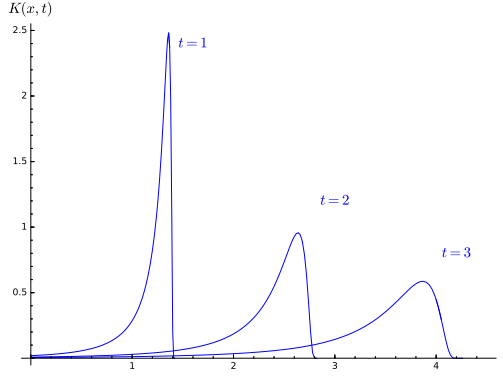}}
 \caption{The wave packet $K(x,t)$, $x\in [0,4.5]$, $t\in \{1, 2, 3\}$.}
 \label{fig-1}
 \end{figure}

Let us write $K_{+}(x,t) = H(x)K(x,t)$. The wave packet $K_+$ for the parameter values $\alpha=\tau=1/2$ is plotted at time instances $t=1, 2, 3$ in Figure \ref{fig-1}. We will interpret $K_{+}$ as a forward moving wave packet with speed $1$. To see why this is justified, consider the rescaled version 
\[
	\K_{t}(\lambda) \coloneqq t K_{+}(\lambda t, t).
\]
For each $t>0$, $\K_{t}$ is a function supported in $[0, 1/\sqrt{\tau}]$ with integral $1/2$.
\begin{proposition}
The function $\K_{t}(\lambda)$ converges to $(1/2)\delta(\lambda-1)$ as $t\to\infty$ in the strong topology of $\MS'(\R)$.
\end{proposition} 
The proof consists essentially of justifying the following heuristic calculation. Sweeping technicalities such as exchanging limits, order of integration, etc.\ under the carpet, we get
\begin{align*}
	\lim_{t\to\infty} \hat{\K}_{t}(\xi) 	
		&= \lim_{t\to\infty} \frac{t}{4\pi\I}\int_{0}^{\infty}\e^{-\I\xi\lambda}\int_{a-\I\infty}^{a+\I\infty}l_{\alpha}(s)\exp\bigl(-\lambda tsl_{\alpha}(s) + st\bigr)\dif s\dif \lambda \\
		&=\lim_{t\to\infty} \frac{t}{4\pi\I}\int_{a-\I\infty}^{a+\I\infty}l_{\alpha}(s)\frac{\e^{ts}}{tsl_{\alpha}(s)+\I\xi}\dif s = \lim_{t\to\infty}\frac{1}{4\pi\I}\int_{a'-\I\infty}^{a'+\I\infty}l_{\alpha}(s/t)\frac{\e^{s}}{sl_{\alpha}(s/t)+\I\xi}\dif s \\
		&=\frac{1}{4\pi\I}\int_{a'-\I\infty}^{a'+\I\infty}\frac{\e^{s}}{s+\I\xi}\dif s = \frac{1}{2}\e^{-\I\xi} = \frac{1}{2}\F\{\delta(\lambda-1)\}(\xi).				
\end{align*}
\begin{proof}
We will show that $\hat{\K}_{t}(\xi) \to (1/2)\e^{-\I\xi}$, boundedly and locally uniformly in $\xi$. Then also $\hat{\K}_{t}(\xi) \to (1/2)\e^{-\I\xi}$ in the strong topology of $\MS'(\R)$. In order to avoid convergency issues of \eqref{eq: K} when $\lambda$ is close to $0$, we first show that for some $\lambda_{0}>0$, $\int_{0}^{\lambda_{0}}\abs{\K_{t}(\lambda)}\dif\lambda \to 0$ as $t\to\infty$. Set
\begin{equation}
\label{eq: L}
L\coloneqq\sup_{\abs{\arg s}\le \pi} \abs{l_{\alpha}(s)}, \quad \lambda_{0} \coloneqq \frac{1}{2L}.
\end{equation}
Note that $L\ge1$, $\lambda_{0}\le1/2$. Using \eqref{eq: S in cone} and changing variables, we have
\[
	\K_{t}(\lambda) = \frac{-1}{4\pi\I}\int_{0}^{\infty}\Bigl\{l_{\alpha}\Bigl(\frac{q}{t}\e^{\I\pi}\Bigr)\e^{\lambda l_{\alpha}(q\e^{\I\pi}/t)q}
											-l_{\alpha}\Bigl(\frac{q}{t}\e^{-\I\pi}\Bigr)\e^{\lambda l_{\alpha}(q\e^{-\I\pi}/t)q}\Bigr\}\e^{-q}\dif q.
\]
For $\lambda\le\lambda_{0}$, $\K_{t}(\lambda)$ is bounded by $\frac{2L}{4\pi}\int_{0}^{\infty}\e^{(1/2)q-q}\dif q = L/\pi$, and converges pointwise to $0$ as $t\to\infty$, since the above integrand is dominated by the integrable function $\e^{-(1/2)q}$, and converges pointwise to $0$. By bounded convergence, it then follows that $\int_{0}^{\lambda_{0}}\K_{t}(\lambda)\e^{-\I\xi\lambda}\dif\lambda$ converges to $0$, uniformly in $\xi$. 

To prove the claim, it then suffices to show that $\int_{\lambda_{0}}^{\infty}\K_{t}(\lambda)\e^{-\I\xi\lambda}\dif\lambda$ converges boundedly and locally uniformly to $(1/2)\e^{-\I\xi}$. Suppose $\xi>0$ (the case $\xi=0$ follows from Proposition \ref{prop: basic properties K}). We use representation \eqref{eq: K} with some $a>0$. Since 
\begin{multline*}
	\int_{a-\I\infty}^{a+\I\infty}\int_{\lambda_{0}}^{\infty}\abs{l_{\alpha}(s)\exp\bigl(-\lambda tsl_{\alpha}(s) + ts -\I\xi\lambda\bigr)}\dif \lambda\, \abs{\dif s} \\
	=\int_{a-\I\infty}^{a+\I\infty}\abs{l_{\alpha}(s)}\frac{\exp\bigl(-\lambda_{0} t\Re(sl_{\alpha}(s)) + at\bigr)}{t\Re(sl_{\alpha}(s))} \abs{\dif s} < \infty,
\end{multline*}
we may interchange the order of integration by the Fubini-Tonelli theorem. We get
\[
	\int_{\lambda_{0}}^{\infty}\K_{t}(\lambda)\e^{-\I\xi\lambda}\dif\lambda 
	= \frac{1}{4\pi\I}\int_{a-\I\infty}^{a+\I\infty}\frac{l_{\alpha}(s)}{sl_{\alpha}(s) + \I\xi/t}\exp\bigl(-\lambda_{0}tsl_{\alpha}(s)+ts - \I\xi\lambda_{0}\bigr)\dif s.
\]
We evaluate this integral by shifting the contour to a Hankel contour encircling the negative real axis, as was done to obtain \eqref{eq: S in cone}. The integral has one singularity, namely the unique zero $s(\xi,t)$ of the function $sl_{\alpha}(s) + \I\xi/t$. Indeed, by applying the argument principle, one sees that this function has a unique zero in the set bounded by the line $[\e^{-\I\pi}R, R]$ and the semicircle with center $0$ and radius $R$ in the lower half plane, and no zeros in the set bounded by the line $[\e^{\I\pi}R, R]$ and the semicircle with center $0$ and radius $R$ in the upper half plane, provided that $\xi>0$ and $R$ is sufficiently large. Since 
\[
	\Re s(\xi,t)l_{\alpha}(s(\xi,t)) = \Re s(\xi, t)\Re l_{\alpha}(s(\xi,t)) - \Im s(\xi,t)\Im l_{\alpha}(s(\xi, t)) = 0,
\]
$\Re l_{\alpha}(s) > 0$ (see \eqref{eq: sgn l_{alpha}}), and $\Im s\Im l_{\alpha}(s) \le 0$, we have that $\Re s(\xi, t) \le 0$.
	
For large $t$, this zero satisfies the asymptotic $s(\xi,t) \sim -\I\xi/t$, as $t\to\infty$, locally uniformly in $\xi$. This can be seen by applying Rouch\'e's theorem on the circle $\abs{s+\I\xi/t} = \abs{\xi}/t^{1+\alpha/2}$. By \eqref{eq: approx l_{alpha} 0}, we have for sufficiently large $t$ that on this circle
\[
	\abs{sl_{\alpha}(s) + \I\xi/t - (s+\I\xi/t)} \lesssim \frac{\abs{\xi}^{1+\alpha}}{t^{1+\alpha}} < \frac{\abs{\xi}}{t^{1+\alpha/2}} = \abs{s+\I\xi t}.
\]
Hence, $sl_{\alpha}(s)+\I\xi/t$ has the same number of zeros (i.e.\ 1) as $s+\I\xi/t$ inside this circle. We get the following representation for sufficiently large $t$:
\begin{align*}
	&\int_{\lambda_{0}}^{\infty}\K_{t}(\lambda)	\e^{-\I\xi\lambda}\dif\lambda \\	
	&= \frac{1}{2}l_{\alpha}(s(\xi, t))\exp\bigl(-\lambda_{0}ts(\xi, t)l_{\alpha}(s(\xi,t)) + ts(\xi, t)-\I\xi\lambda_{0}\bigr) \\
	&\quad +\frac{\e^{-\I\xi\lambda_{0}}}{4\pi\I}\int_{0}^{\infty}\biggl\{\frac{l_{\alpha}(q\e^{\I\pi})\e^{\lambda_{0}tql_{\alpha}(q\e^{\I\pi})}}{\I\xi/t - ql_{\alpha}(q\e^{\I\pi})} 
												- \frac{l_{\alpha}(q\e^{-\I\pi})\e^{\lambda_{0}tql_{\alpha}(q\e^{-\I\pi})}}{\I\xi/t - ql_{\alpha}(q\e^{-\I\pi})}\biggr\}\e^{-tq}\dif q\\
	&= \frac{1}{2}l_{\alpha}(s(\xi,t))\e^{ts(\xi,t)}
			 + \frac{\e^{-\I\xi\lambda_{0}}}{4\pi\I}\int_{0}^{\infty}\biggl\{\frac{l_{\alpha}(\e^{\I\pi}q/t)\e^{\lambda_{0}ql_{\alpha}(\e^{\I\pi}q/t)}}{\I\xi - ql_{\alpha}(\e^{\I\pi}q/t)} - 
			 \frac{l_{\alpha}(\e^{-\I\pi}q/t)\e^{\lambda_{0}ql_{\alpha}(\e^{-\I\pi}q/t)}}{\I\xi - ql_{\alpha}(\e^{-\I\pi}q/t)}\biggr\}\e^{-q}\dif q.										
\end{align*}
Since $\Re s(\xi,t)\le0$, the first term is bounded, and it converges locally uniformly to $(1/2)\e^{-\I\xi}$, in view of the asymptotic $s(\xi,t)\sim-\I\xi/t$. The integral converges uniformly to $0$ as $t\to\infty$. Given $\eps>0$, one can first find some $Q$ so that 
\[
	\int_{Q}^{\infty}2\frac{L\e^{\lambda_{0}Lq}}{(\sqrt{\tau}/2) q}\e^{-q}\dif q \le  \frac{4L}{\sqrt{\tau}}\int_{Q}^{\infty}\frac{\e^{-(1/2)q}}{q}\dif q \le \frac{\eps}{2}.
\]
On the interval $[0,Q]$, we will use a Taylor approximation. We get
\[
	\frac{l_{\alpha}(\e^{\pm\I\pi}q/t)\e^{\lambda_{0}ql_{\alpha}(\e^{\pm\I\pi}q/t)}}{\I\xi-ql_{\alpha}(\e^{\pm\I\pi}q/t)} = \frac{\e^{\lambda_{0}q}}{\I\xi-q}
	\biggl\{1+O_{Q}\biggl(\frac{q^{\alpha}}{t^{\alpha}} + \frac{q^{1+\alpha}}{t^{\alpha}} + \frac{q^{1+\alpha}}{t^{\alpha}\abs{\I\xi-q}}\biggr)\biggr\}.	
\]
Let then $t$ be so large that 
\[
	\abs{\frac{l_{\alpha}(\e^{\I\pi}q/t)\e^{\lambda_{0}ql_{\alpha}(\e^{\I\pi}q/t)}}{\I\xi-ql_{\alpha}(\e^{\I\pi}q/t)} - \frac{l_{\alpha}(\e^{-\I\pi}q/t)\e^{\lambda_{0}ql_{\alpha}(\e^{-\I\pi}q/t)}}{\I\xi-ql_{\alpha}(\e^{-\I\pi}q/t)}} \le \frac{\alpha \eps}{2Q^{\alpha}} q^{\alpha-1}, \quad \text{ for } q\in [0,Q].
\]
Then $\abs[1]{\int_{0}^{Q} \dotso}\le \eps/2$. We conclude that $\int_{\lambda_{0}}^{\infty}\K_{t}(\lambda)	\e^{-\I\xi\lambda}\dif\lambda$ converges boundedly and locally uniformly to $(1/2)\e^{-\I\xi}$, which finishes the proof of the proposition.
\end{proof}

This proposition gives some indication that $K_{+}$ is concentrated around $x=t$. This ``concentration'' around $x=t$ is however much less drastic than the concentration of $\K_{t}$ around $\lambda=1$. It is for example not the case that $K_{+}(x,t) - (1/2)\delta(x-t) \to 0$. Actually, the wave packet will spread out in space, albeit on a scale smaller than $\abs{x-t} \simeq t$. Namely, we will see that $K_{+}(x,t)$ can be described as a wave packet of height $\simeq t^{-\frac{1}{1+\alpha}}$ and width $\simeq t^{\frac{1}{1+\alpha}}$ centered around $x=t$. Let us first give a dispersion estimate for $K_{+}(x,t)$. For later use, we also bound the derivatives with respect to $x$.

\begin{proposition}
\label{prop: uniform bound K}
For every $n\in \N$, we have the bound
\[
	\norm{\dpd[n]{K_{+}}{x}}_{L^{\infty}_{x}} \lesssim_{n} t^{-\frac{n+1}{1+\alpha}}.
\]
\end{proposition}
\begin{proof}
Let $L$ and $\lambda_{0}$ be as before, see \eqref{eq: L}. Suppose first that $0\le x \le \lambda_{0}t$. Then
\begin{align*}
	\dpd[n]{K_{+}}{x}(x,t) &= \frac{1}{4\pi\I}\int_{0}^{\infty}q^{n}\Bigl(l_{\alpha}(q\e^{-\I\pi})^{n+1}\e^{xql_{\alpha}(q\e^{-\I\pi})} -l_{\alpha}(q\e^{\I\pi})^{n+1}\e^{xql_{\alpha}(q\e^{\I\pi})}\Bigr)\e^{-qt}\dif q\\
					&\lesssim_{n}\int_{0}^{\infty}q^{n}\e^{xLq-tq}\dif q \le \int_{0}^{\infty}q^{n}\e^{-(1/2)tq}\dif q \lesssim_{n} t^{-n-1}.
\end{align*}
For $x\ge \lambda_{0}t$, we use representation \eqref{eq: K} and move the contour to the imaginary axis. We get
\[
	\dpd[n]{K_{+}}{x}(x,t)	= \frac{1}{4\pi}\int_{-\infty}^{\infty}(-\I y)^{n}l_{\alpha}(\I y)^{n+1}\exp\bigl(-x\I yl_{\alpha}(\I y) + t\I y\bigr)\dif y 
					\lesssim_{n} \int_{0}^{\infty}y^{n}\exp\bigl(xy\Im l_{\alpha}(\I y)\bigr)\dif y.
\]
By \eqref{eq: bound l_{alpha}}, we get
\[
	\dpd[n]{K_{+}}{x}(x,t)	\lesssim \int_{0}^{1}y^{n}\exp\bigl(-\lambda_{0}tc_{1}y^{1+\alpha}\bigr)\dif y + \int_{1}^{\infty}y^{n}\exp\bigl(-\lambda_{0}tc_{2}y^{1-\alpha}\bigr)\dif y 
					\lesssim_{n} t^{-\frac{n+1}{1+\alpha}}.
\]
\end{proof}

\begin{figure}[htb]
   \begin{center}$
   \begin{array}{ccc}
   \includegraphics[width=0.45\textwidth]{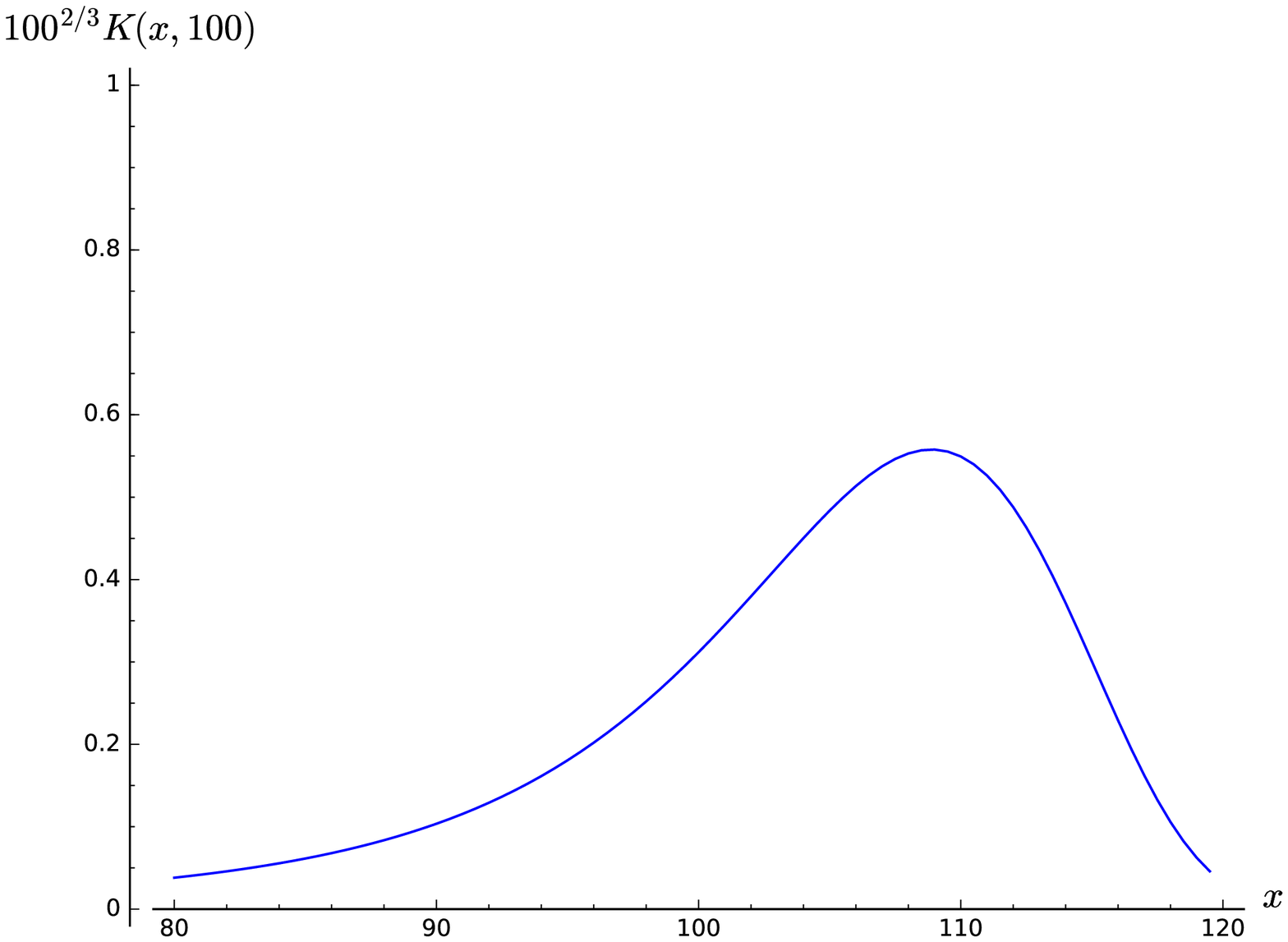}
   \hphantom{blabla}
   \includegraphics[width=0.45\textwidth]{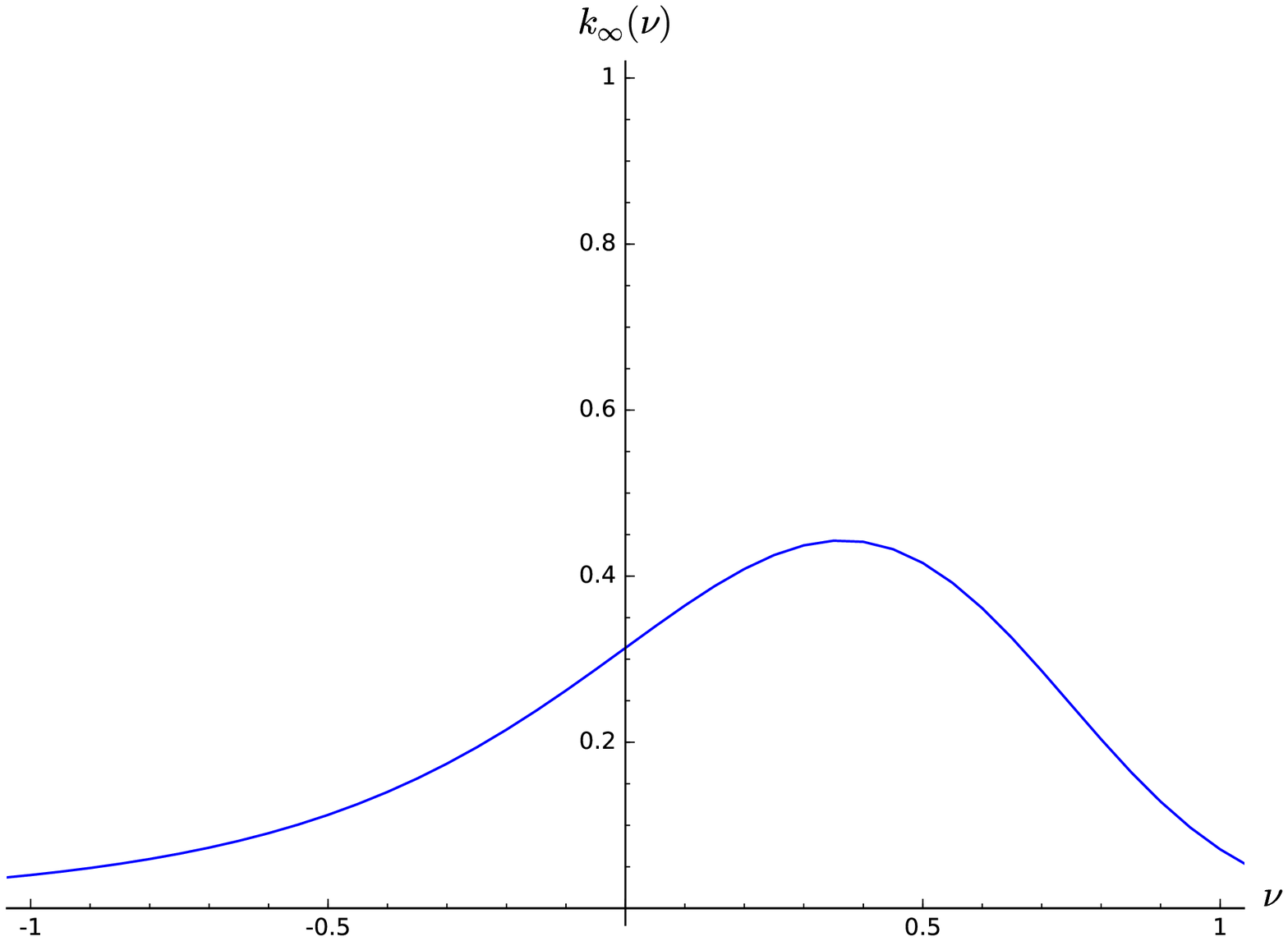} 
   \end{array}$
   \end{center}
   \caption{Comparison between $100^{\frac{1}{1+\alpha}}K(x,100)$ and $k_\infty (\nu)$. } \label{figure}
   \end{figure}

\begin{figure}[htb]
 \centerline{\includegraphics[width=7cm]{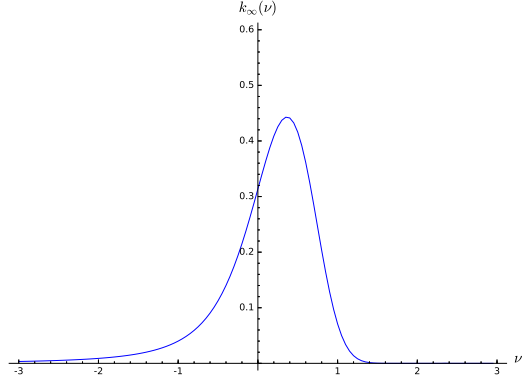}}
 \caption{The function $k_\infty$ for $\nu\in [-3,3]$.}
 \label{fig-3}
 \end{figure}

We will now give a precise description of the shape of the wave packet, in the limit $t\to\infty$. For this, we introduce the function
\[
	k_{t}(\nu) \coloneqq t^{\frac{1}{1+\alpha}}K_{+}(t+\nu t^{\frac{1}{1+\alpha}}, t), \quad \nu\in \R.
\]
\begin{theorem}
\label{thm: shape wave packet}
There exists a function $k_{\infty}(\nu)$ with the property that $k_{t}(\nu) \to k_{\infty}(\nu)$ as $t\to\infty$, locally uniformly in $\nu$.
The function  $k_{\infty}$ has the following representation:
\begin{equation}
\label{eq: k_{infty}}
	k_{\infty}(\nu) = \frac{1}{4\pi}\int_{-\infty}^{\infty}\exp\Bigl(\frac{1-\tau}{2}(\I w)^{1+\alpha} - \I \nu w\Bigr)\dif w.
\end{equation}
Also, $\partial_{\nu}^{n}k_{t}(\nu) \to k^{(n)}_{\infty}(\nu)$ locally uniformly in $\nu$, for every $n\in \N$.
\end{theorem}
\begin{proof}
We again use representation \eqref{eq: K} on the imaginary axis. We split the range of integration into three parts as follows: 
\begin{align*}
	k_{t}(\nu) 	&= \frac{t^{\frac{1}{1+\alpha}}}{4\pi\I}\int_{-\I\infty}^{\I\infty}l_{\alpha}(s)\exp\bigl(ts(1-l_{\alpha}(s)) - \nu t^{\frac{1}{1+\alpha}}sl_{\alpha}(s)\bigr)\dif s \\
			&= \frac{t^{\frac{1}{1+\alpha}}}{4\pi} \biggl( \int_{\abs{y}\le t^{-\frac{1}{1+2\alpha}}} + \int_{t^{-\frac{1}{1+2\alpha}} \le\, \abs{y} \le 1} + \int_{\abs{y }\ge 1} \biggr)
						l_{\alpha}(\I y)\exp\bigl(t\I y(1-l_{\alpha}(\I y)) - \nu t^{\frac{1}{1+\alpha}}\I y l_{\alpha}(\I y)\bigr)\dif y \\
			&\eqqcolon I_{1} + I_{2} + I_{3}.
\end{align*}

When $0\le y \le 1$, we have $\Im l_{\alpha}(\I y) \le -c_{1}y^{\alpha}$, see \eqref{eq: bound l_{alpha}}. By \eqref{eq: approx l_{alpha} 0}, we also have $\abs{\Im l_{\alpha}(\I y)} \le \tilde{c}_{1}y^{\alpha}$ for some $\tilde{c}_{1}>0$. This yields
\begin{align*}
	I_{2} 	&\lesssim t^{\frac{1}{1+\alpha}}\int_{t^{-\frac{1}{1+2\alpha}}}^{1}\exp\bigl(-tc_{1}y^{1+\alpha} + \abs{\nu}t^{\frac{1}{1+\alpha}}\tilde{c}_{1}y^{1+\alpha}\bigr)\dif y \\
		&\lesssim t^{\frac{1}{1+\alpha}}\int_{t^{-\frac{1}{1+2\alpha}}}^{1}\exp\bigl(-t(c_{1}/2)y^{1+\alpha}\bigr) \dif y 
		\lesssim \int_{(c_{1}/2)^{\frac{1}{1+\alpha}}t^{\frac{\alpha}{(1+\alpha)(1+2\alpha)}}}^{(tc_{1}/2)^{\frac{1}{1+\alpha}}}\e^{-w^{1+\alpha}}\dif w \to 0,
\end{align*}
as $t\to\infty$. Here we used that $\tilde{c}_{1}\abs{\nu}t^{\frac{1}{1+\alpha}} \le c_{1}t/2$ for sufficiently large $t$, uniformly for $\nu$ in compact sets.

When $y\ge1$, we have $\Im l_{\alpha}(\I y) \le -c_{2}y^{-\alpha}$ and $\abs{\Im l_{\alpha}(\I y)} \le \tilde{c}_{2}y^{-\alpha}$ for some $\tilde{c}_{2}>0$ (see again \eqref{eq: bound l_{alpha}} and \eqref{eq: approx l_{alpha}}). This implies that 
\begin{align*}
	I_{3}	&\lesssim t^{\frac{1}{1+\alpha}}\int_{1}^{\infty}\exp\bigl(-tc_{2}y^{1-\alpha} + \abs{\nu}t^{\frac{1}{1+\alpha}}\tilde{c}_{2}y^{1-\alpha}\bigr)\dif y \\
		&\lesssim t^{\frac{1}{1+\alpha}}\int_{1}^{\infty}\exp\bigl(-t(c_{2}/2)y^{1-\alpha}\bigr)\dif y 
		\lesssim t^{\frac{1}{1+\alpha}-\frac{1}{1-\alpha}}\int_{(tc_{2}/2)^{\frac{1}{1-\alpha}}}^{\infty}\e^{-w^{1-\alpha}}\dif w \to 0,
\end{align*}
as $t\to\infty$, uniformly for $\nu$ in compact sets. Hence, in the limit $t\to\infty$, only the contribution from $I_{1}$ remains.

To treat $I_{1}$, we will approximate $l_{\alpha}$ using \eqref{eq: approx l_{alpha} 0}. For the remainders, we write 
\[
	\tilde{E}_{1}(s) = l_{\alpha}(s) - 1, \quad \tilde{E}_{2}(s) \coloneqq l_{\alpha}(s) - \Bigl(1 - \frac{1-\tau}{2}s^{\alpha}\Bigr).
\]
In the integral $I_{1}$, change variables to $w = t^{\frac{1}{1+\alpha}}y$ and approximate $l_{\alpha}(\I y)$ by $1 - (1/2)(1-\tau)(\I y)^{\alpha}$. This gives
\[
	I_{1} = \frac{1}{4\pi}\int_{-t^{\beta}}^{t^{\beta}}
		l_{\alpha}(t^{-\frac{1}{1+\alpha}}\I w)\exp\Bigl(\frac{1-\tau}{2}(\I w)^{1+\alpha} - \I \nu w 
		- t^{\frac{\alpha}{1+\alpha}}\I w\tilde{E}_{2}(t^{-\frac{1}{1+\alpha}}\I w) - \nu\I w\tilde{E}_{1}(t^{-\frac{1}{1+\alpha}}\I w)\Bigr)\dif w.
\]
Here, $\beta=\frac{\alpha}{(1+\alpha)(1+2\alpha)}$. Since $\tilde{E}_{1}(s) \lesssim \abs{s}^{\alpha}$ and $\tilde{E}_{2}(s) \lesssim \abs{s}^{2\alpha}$ for $s\to0$, the integrand converges pointwise to the function ${\exp\bigl((1/2)(1-\tau)(\I w)^{1+\alpha} - \I\nu w\bigr)}$, uniformly for $\nu$ in compact sets. Furthermore, on the interval $[-t^{\beta}, t^{\beta}]$, 
\[
	\nu\I w\tilde{E}_{1}(t^{-\frac{1}{1+\alpha}}\I w) \lesssim t^{-\frac{\alpha^{2}}{(1+\alpha)(1+2\alpha)}} \lesssim 1, \quad t^{\frac{\alpha}{1+\alpha}}\I w\tilde{E}_{2}(t^{-\frac{1}{1+\alpha}}\I w) \lesssim 1,
\]
so the integrand is dominated by the integrable function $\exp\bigl(-(1/2)(1-\tau)\sin(\alpha\pi/2)\abs{w}^{1+\alpha}\bigr)$. By dominated convergence,
\[
	I_{1}\to \frac{1}{4\pi}\int_{-\infty}^{\infty}\exp\Bigl(\frac{1-\tau}{2}(\I w)^{1+\alpha} - \I \nu w\Bigr)\dif w \eqqcolon k_{\infty}(\nu),
\]
as $t\to\infty$, uniformly for $\nu$ in compact sets.

The proof for $\partial^{n}_{\nu}k_{t}(\nu)$ is completely analogous. The corresponding integrals $I_{2}$ and $I_{3}$ tend to zero, while the integral $I_{1}$ converges to 
\[
	\frac{1}{4\pi}\int_{-\infty}^{\infty}(-\I w)^{n}\exp\Bigl(\frac{1-\tau}{2}(\I w)^{1+\alpha} - \I \nu w\Bigr)\dif w = \dod[n]{}{\nu}k_{\infty}(\nu).
\]
\end{proof}
In Figure \ref{figure} we compare the shape of the wave packet at some large time with the function $k_\infty (\nu)$. We choose again parameter values $\alpha=\tau=1/2$. On the left we show $K$ at time $t=100$ scaled by factor $t^{\frac{1}{1+\alpha} }= 100^{\frac{2}{3}}\approx 21.54$,  for $x\in [80,120]$. On the right we show a plot of the function $k_{\infty}$ for $\nu\in [-1,1]$. In Figure  \ref{fig-3} we plot  $k_{\infty}$ in a larger range.
 
Let us list some properties of the function $k_{\infty}(\nu)$. 
\begin{itemize}
	\item From the representation \eqref{eq: k_{infty}}, it follows immediately that $k_{\infty}$ belongs to the Gevrey class $G^{\frac{1}{1+\alpha}}(\R)$. In particular, it is an entire function. 
	\item The function $k_{\infty}$ is real valued, since it can be written as 
\[
	k_{\infty}(\nu) = \frac{1}{2\pi}\Re \int_{0}^{\infty}\exp\Bigl(\frac{1-\tau}{2}(\I w)^{1+\alpha} - \I \nu w\Bigr)\dif w.
\]
The integral from $0$ to $\infty$ above also defines an entire function. It belongs to a general class of Fourier-Laplace transforms which was investigated in \cite{BDV2021}. In that paper, a complete asymptotic analysis of such functions was performed on half-lines emanating from the origin. Taking the first term from the asymptotic series \cite[Equation (4.1)]{BDV2021} with non-vanishing real part, we see that 
\begin{equation}
\label{eq: asymptotics k-}
	k_{\infty}(\nu) \sim \frac{\sin(\alpha\pi)}{2\pi}\Bigl(\frac{1-\tau}{2}\Bigr)^{\frac{1}{1+\alpha}}\frac{\Gamma(2+\alpha)}{\abs{\nu}^{2+\alpha}}, \quad \text{as } \nu \to -\infty.
\end{equation}
For $\nu \to +\infty$, every term in the asymptotic series \cite[Equation (4.1)]{BDV2021} is purely imaginary, so this only tells us that $k_{\infty}(\nu) \lesssim_{n} \nu^{-n}$, as $\nu \to \infty$, for every $n$. However, using the saddle-point method, one can determine its precise asymptotic. We restrict ourselves here to just sketching the method. Suppose $\nu>0$. In \eqref{eq: k_{infty}}, change variables $w = \nu^{\frac{1}{\alpha}}y$ to get 
\[
	k_{\infty}(\nu) = \frac{\nu^{\frac{1}{\alpha}}}{4\pi}\int_{-\infty}^{\infty}\exp\Bigl\{\nu^{1+\frac{1}{\alpha}}\Bigl( \frac{1-\tau}{2}(\I y)^{1+\alpha} - \I y\Bigr)\Bigr\}\dif y.
\]
The function $((1-\tau)/2) z^{1+\alpha} - z$ is holomorphic on $\C\setminus (-\infty,0]$ and has a unique saddle point $z_{0} = \bigl((1-\tau)(1+\alpha)/2\bigr)^{-\frac{1}{\alpha}}$. One can shift the contour of integration to a contour passing through this saddle point via the ``steepest path.'' Applying the saddle point method then yields the following asymptotic:
\begin{equation}
\label{eq: asymptotics k+}
	k_{\infty}(\nu) \sim \frac{1}{4}\sqrt{\frac{2z_{0}}{\alpha\pi}}\nu^{\frac{1}{2\alpha}-\frac{1}{2}}
		\exp\Bigl(-\frac{\alpha z_{0}}{1+\alpha}\nu^{1+\frac{1}{\alpha}}\Bigr), \quad \text{as } \nu \to \infty.
\end{equation}
	\item Another interesting property is that $k_{\infty}$ is a close cousin of the Gaussian and Airy functions. Indeed, renormalizing by setting 
\[
	k(\nu)\coloneqq \biggl(\frac{(1-\tau)(1+\alpha)}{2}\biggr)^{\frac{1}{1+\alpha}}k_{\infty}\biggl(-\biggl(\frac{(1-\tau)(1+\alpha)}{2}\biggr)^{\frac{1}{1+\alpha}}\nu\biggr),
\]
we have that $k$ satisfies the fractional ordinary differential equation
\[
	\fracD{-\infty}{\alpha}{\nu}k(\nu) + \nu k(\nu) = 0.
\]
This follows immediately by taking Fourier transforms, since $\hat{k}(\xi) = \exp\bigl(\frac{1}{1+\alpha}(\I\xi)^{1+\alpha} \bigr)$. 
	\item From the previous property we can deduce that $k$ and hence also $k_{\infty}$ is everywhere positive. From both asymptotics \eqref{eq: asymptotics k-} and \eqref{eq: asymptotics k+} we see that $k$ is eventually positive and so it has (at most) finitely many zeros. Suppose that it does have zeros. Let $\nu_{0}$ be the smallest one. Then $k(\nu)>0$ for every $\nu<\nu_{0}$. From the above differential equation, it follows that 
\[
	\fracD{-\infty}{\alpha}{\nu}k(\nu_{0}) = \frac{1}{\Gamma(-\alpha)}\Fp \int_{-\infty}^{\nu_{0}}k(\nu)(\nu_{0}-\nu)^{-\alpha-1}\dif \nu = 0.
\]
Here we used the Hadamard finite part to compute the $\alpha$-th order fractional derivative. However, if $\nu_{0}$ is a zero of $k$, then $k(\nu)(\nu_{0}-\nu)^{-\alpha-1}$ is integrable. Hence we get
\[
	0 = \int_{-\infty}^{\nu_{0}} k(\nu)(\nu_{0}-\nu)^{-\alpha-1}\dif \nu > 0,
\]
a contradiction. We conclude that $k$ has no zeros, so it is everywhere positive.
\end{itemize}

\begin{remark}
Theorem \ref{thm: shape wave packet} can be rephrased by saying that $K_{+}(t+\nu t^{\frac{1}{1+\alpha}}, t) \sim k_{\infty}(\nu)t^{-\frac{1}{1+\alpha}}$, as $t\to\infty$, locally uniformly in $\nu$. In particular we have 
\[
	K_{+}(t,t) \sim k_{\infty}(0)t^{-\frac{1}{1+\alpha}} = \frac{1}{2\pi(1+\alpha)}\sin\Bigl(\frac{\pi}{1+\alpha}\Bigr)\Gamma\Bigl(\frac{1}{1+\alpha}\Bigr)\Bigl(\frac{2}{1-\tau}\Bigr)^{\frac{1}{1+\alpha}}
		t^{-\frac{1}{1+\alpha}}.
\]
It is also possible to determine the asymptotics of $K_{+}(\lambda t, t)$ with $\lambda\neq1$. If $0\le\lambda<1$, we have the power decay
\[
	K_{+}(\lambda t, t) \sim \frac{\sin(\alpha\pi)(1-\tau)(1+2\lambda)}{4\pi(1-\lambda)^{2+\alpha}}\Gamma(1+\alpha)t^{-1-\alpha}.
\]
This asymptotic relation holds uniformly for $\lambda\in[0,\lambda_{1}]$, for any $\lambda_{1}<1$. If $1<\lambda<1/\sqrt{\tau}$, we have exponential decay. Set $f(s) = s-\lambda sl_{\alpha}(s)$, and denote by $a_{\lambda}$ the unique positive real zero of the function $f'(s)$. Then
\[
	K_{+}(\lambda t, t) \sim \frac{l_{\alpha}(a_{\lambda})}{4}\sqrt{\frac{2}{\pi f''(a_{\lambda})}}\frac{\e^{f(a_{\lambda})t}}{\sqrt{t}}.
\]
We remark that $f(a_{\lambda}) < 0$ and $f''(a_{\lambda})>0$, and that $f(a_{\lambda}) \to -\infty$ if $\lambda\to1/\sqrt{\tau}$.

Both of these asymptotic relations can be obtained via the method of steepest descent, but we omit the details.
\end{remark}

In view of the preceding discussion, it is natural to consider $K_{+}$ as a dispersive wave packet with speed $1$ and wave front speed $1/\sqrt{\tau}$. In previous works on (fractional) wave equations, several ways of assigning a velocity to waves in dissipative media are used. See for example \cite{Luchko2013}. If one defines the maximum position, the center of gravity, and the center of mass of the wave respectively as
\begin{align*}
	x^{\textrm{max}}(t) 	&= \argmax_{x}K_{+}(x,t); \\
	x^{g}(t) 			&= \frac{\int_{0}^{\infty}xK_{+}(x,t)\dif x}{\int_{0}^{\infty}K_{+}(x,t)\dif x}; \\
	x^{m}(t) 			&= \frac{\int_{0}^{\infty}xK^{2}_{+}(x,t)\dif x}{\int_{0}^{\infty}K^{2}_{+}(x,t)\dif x};
\end{align*}
one can define associated velocities as the instantaneous or average propagation speed of these points. In our case, it would appear that
\[
	x^{\textrm{max}}(t)\sim t, \quad x^{g}(t)\sim t, \quad x^{m}(t) \sim t, 
\]
so the associated velocities would all be  asymptotically equal to $1$. 

\begin{remark} 
Note that \eqref{eq: FZWE} resulted from a reduction to dimensionless quantities (i.e.\ \eqref{eq:system}). For the model in its original form, so with Newton's second law in the form $\partial_{x}\sigma = \rho\partial_{t}^{2}u$, including the density constant $\rho$, and with the fractional Zener constitutive law in the form \eqref{eq:frac Zener model}, the wave front speed and wave packet speed are given by $\sqrt{E\tau_{\eps}/(\rho\tau_{\sigma})}$ and $\sqrt{E/\rho}$ respectively. These speeds can be related to the limiting values of the material functions of the body. We have
\begin{align*}
	\text{wave front speed } 	&= \sqrt{\frac{E\tau_{\eps}}{\rho\tau_{\sigma}}}	= \frac{1}{\sqrt{\rho J_{g}}} = \sqrt{\frac{G_{g}}{\rho}}, \\
	\text{wave packet speed } &= \sqrt{\frac{E}{\rho}} = \frac{1}{\sqrt{\rho J_{e}}} = \sqrt{\frac{G_{e}}{\rho}}.	 
\end{align*} 
Here, $J_{g}$ and $G_{g}$ are the glass compliance and glass modulus, related to the instantaneous behavior of the material, and $J_{e}$ and $G_{e}$ are the equilibrium compliance and equilibrium modulus, related to the equilibrium behavior of the material, see e.g.\ \cite[Chapter 2]{Mainardi}.
In dimensionless form  $G_{g}= 1/ J_{g} = 1/\tau$ and $G_{e}=1/J_{e}=1$.  For a more general class of materials (including the fractional Zener model) they are calculated  and presented in \cite[Table 1]{ZO2020}. 
\end{remark}
	
Finally, let us describe the shape of the solution with initial conditions $u(x,0)=u_{0}(x)$, $\partial_{t}u(x,0)=0$, given by $u(x,t) = K(x,t)\ast_{x}u_{0}(x)$.
\begin{theorem}
Suppose $u_{0}\in \MS$ and suppose that $\int u_{0}(x)\dif x \neq 0$. Then
\[
	\norm{u}_{L^{\infty}_{x}} \lesssim t^{-\frac{1}{1+\alpha}}\norm[1]{u_{0}}_{L^{1}_{x}},
\]
and $u(x,t) = K_{+}(x,t)\ast_{x}u_{0}(x) + K_{+}(-x,t)\ast_{x}u_{0}(x) \eqqcolon u_{+}(x,t) + u_{-}(x,t)$, where
\[
	t^{\frac{1}{1+\alpha}}u_{+}(t+\nu t^{\frac{1}{1+\alpha}}, t) \to Ak_{\infty}(\nu), \quad t^{\frac{1}{1+\alpha}}u_{-}(-t-\nu t^{\frac{1}{1+\alpha}},t) \to Ak_{\infty}(\nu), \quad A = \int_{-\infty}^{\infty}u_{0}(x)\dif x,
\]
as $t\to\infty$, locally uniformly in $\nu$.
\end{theorem}
\begin{proof}
By Proposition \ref{prop: uniform bound K},
\[
	u(x,t) = \int_{x-t/\sqrt{\tau}}^{x+t/\sqrt{\tau}}K(x-y,t)u_{0}(y) \dif y \lesssim t^{-\frac{1}{1+\alpha}}\int_{-\infty}^{\infty}\abs[0]{u_{0}(x)}\dif x.
\]
Set now $x=t+\nu t^{\frac{1}{1+\alpha}}$. Rewriting the convolution in terms of $k_{t}(\nu)$ gives 
\[
	u_{+}(x,t) = \int K_{+}(x-y,t)u_{0}(y)\dif y = t^{-\frac{1}{1+\alpha}}\int k_{t}(\nu-yt^{-\frac{1}{1+\alpha}})u_{0}(y)\dif y.
\]
Applying Theorem \ref{thm: shape wave packet} and dominated convergence, we see that locally uniformly in $\nu$
\[
	t^{\frac{1}{1+\alpha}}u_{+}(t+\nu t^{\frac{1}{1+\alpha}}, t) \to \int k_{\infty}(\nu)u_{0}(y)\dif y = \biggl(\int_{-\infty}^{\infty}u_{0}(y)\dif y\biggr) k_{\infty}(\nu).
\]
The proof for $u_{-}$ is analogous.
\end{proof}
\begin{remark}
If $\int u_{0}(x)\dif x = 0$ and $\int x u_{0}(x)\dif x \neq 0$, than $u_{0}$ has a primitive $u_{0}^{(-1)}$ in $\MS$. We can integrate by parts in the convolution:
\[
	u(x,t) = \int_{x-t/\sqrt{\tau}}^{x+t/\sqrt{\tau}}K(x-y,t)u_{0}(y) \dif y = \int_{x-t/\sqrt{\tau}}^{x+t/\sqrt{\tau}}\dpd{K}{x}(x-y,t)u_{0}^{(-1)}(y)\dif y,
\]
since the boundary terms vanish. Similarly as in the above proof, using Proposition \ref{prop: uniform bound K} and Theorem \ref{thm: shape wave packet}, now defining $u_{+}(x,t) = \partial_{x}K_{+}(x,t)\ast_{x}u_{0}^{(-1)}(x)$ and $u_{-}(x,t) = \partial_{x}K_{+}(-x,t)\ast_{x}u_{0}^{(-1)}(x)$,
\begin{gather*}
	\norm{u}_{L^{\infty}_{x}} \lesssim t^{-\frac{2}{1+\alpha}}\norm[0]{u_{0}^{(-1)}}_{L^{1}_{x}}, \\ 
	 t^{\frac{2}{1+\alpha}}u_{+}(t+\nu t^{\frac{1}{1+\alpha}}, t) \to \biggl(\int_{-\infty}^{\infty} u_{0}^{(-1)}(x)\dif x\biggr) k_{\infty}'(\nu) = \biggl(-\int_{-\infty}^{\infty}xu_{0}(x)\dif x\biggr)k_{\infty}'(\nu),
\end{gather*}
and similarly for $u_{-}$.

If more moments of $u_{0}$ vanish, then integrating by parts introduces extra boundary terms, by the non-differentiability of $K$ at $x=0$. Suppose $n\ge 2$ is the smallest integer such that $\int x^{n}u_{0}(x)\dif x \neq 0$. Denoting the $j$-th order primitive, $j\le n$, of $u_{0}$ in $\MS$ by $u_{0}^{(-j)}$, and setting $m=2\lfloor n/2\rfloor$, we get
\[
	u(x,t) = 2u_{0}^{(-2)}(x)\dpd{K}{x}(0^{+}, t) + \dotsb + 2u_{0}^{(-m)}(x)\dpd[m-1]{K}{x}(0^{+}, t) + \int_{x-t/\sqrt{\tau}}^{x+t/\sqrt{\tau}}\dpd[n]{K}{x}(x-y,t)u_{0}^{(-n)}(y)\dif y.
\]
It is possible to estimate $\partial_{x}^{j}K(0^{+}, t)$ for large $t$; using these estimates, one might then estimate the $L^{\infty}_{x}$-norm of $u$ by a linear combination of the $L^{\infty}_{x}$-norms of $u_{0}^{(-2)}, \dotso, u_{0}^{(-m)}$ and the $L^{1}_{x}$-norm of $u_{0}^{(-n)}$, where the coefficients are negative powers of $t$ depending on the order of the respective primitive of $u_{0}$ and $\alpha$.
\end{remark}
\section{The case $\alpha=1$}\label{ClassZ}
We now briefly discuss the case $\alpha=1$, which is known as the (classical) Zener model, or the Standard Linear Solid (SLS) model. The SLS wave equation is 
\begin{equation}
\label{eq: CZWE}
	\dpd[2]{}{t} u(x,t) = \ML^{-1}_{s\to t} \left(\frac{1+s}{1+\tau s}\right) \ast_{t} \dpd[2]{}{x} u(x,t) ,\qquad x\in \mathbb{R}, \quad t>0,
\end{equation}
so \eqref{eq: FZWE} with $\alpha=1$. The fundamental solution $S$ of \eqref{eq: CZWE} is again supported in the forward cone $\abs{x}\le t/\sqrt{t}$, but we will see that it is not smooth on the boundary of this cone, in contrast with the case $0<\alpha<1$. The Laplace transform $\tilde{S}$ of $S$ is now given by
\[
	\tilde{S}(x,s) = \frac{1}{2s}\sqrt{\frac{1+\tau s}{1+s}}\exp\biggl(-\abs{x} s \sqrt{\frac{1+\tau s}{1+s}}\biggr).
\]
Note that this function has analytic continuation to $\C \setminus\bigl( \{0\} \cup [-1/\tau, -1]\bigr)$. The point $s=0$ is a simple pole of $\tilde{S}$; the line segment $[-1/\tau,-1]$ is a branch cut.
\begin{theorem}
The fundamental solution $S$ of \eqref{eq: CZWE} is discontinuous at the boundary of the forward light cone. More precise, it has the following form:
\begin{align*}
	S(x,t)	&= \frac{\sqrt{\tau}}{2}\exp\biggl(-\frac{\sqrt{\tau}}{2}\Bigl(\frac{1}{\tau}-1\Bigr)\abs{x}\biggr)H(t-\sqrt{\tau}\abs{x}) + E(x,t),
\end{align*}
where $E$ is a continuous function supported in the forward cone $\{(x,t): \abs{x} \le t/\sqrt{\tau}\}$.
\end{theorem}
\begin{proof}
The fact that $S(x,t)$, and hence also $E(x,t) \coloneqq S(x,t) - (\sqrt{\tau}/2)\e^{-\frac{\sqrt{\tau}}{2}(1/\tau-1)\,\abs{x}}{H(t-\sqrt{\tau}\abs{x})}$, is supported in the cone $\abs{x}\le t/\sqrt{\tau}$ can be proven in a similar fashion as in Proposition \ref{prop: support S}. (In fact, that proof does not require that $\alpha<1$.)

To show that $E$ is continuous, we first note that, similarly as in \eqref{eq: approx l_{alpha}}, 
\begin{equation}
\label{eq: approx l_{1}}
	l_{1}(s) \coloneqq \sqrt{\frac{1+\tau s}{1+s}} = \sqrt{\tau}\biggl(1 + \frac{1}{2}\Bigl(\frac{1}{\tau}-1\Bigr)s^{-1} + O(\,\abs{s}^{-2})\biggr), \quad \text{as } \abs{s} \to \infty.
\end{equation}
We have 
\begin{align*}
	\ML\{E(x,t)\}(s) 	&= \frac{1}{2s}l_{1}(s)\exp\bigl(-\abs{x} sl_{1}(s)\bigr) 
				- \frac{\sqrt{\tau}}{2s}\exp\biggl(-\sqrt{\tau}\abs{x} s - \frac{\sqrt{\tau}}{2}\Bigl(\frac{1}{\tau}-1\Bigr)\abs{x}\biggr) \\
				&= \frac{\sqrt{\tau}}{2s}\exp\biggl(-\sqrt{\tau}\abs{x} s - \frac{\sqrt{\tau}}{2}\Bigl(\frac{1}{\tau}-1\Bigr)\abs{x}\biggr)\Bigl(\bigl(1+O(1/\abs{s})\bigr)\exp\bigl(O(\abs[0]{x/s})\bigr) - 1\Bigr) \\
				&= \frac{\sqrt{\tau}}{2s}\exp\biggl(-\sqrt{\tau}\abs{x} s - \frac{\sqrt{\tau}}{2}\Bigl(\frac{1}{\tau}-1\Bigr)\abs{x}\biggr) \cdot O\Bigl(\frac{1+\abs{x}}{\abs{s}}\Bigr), \quad \text{as } \abs{s}\to\infty.
\end{align*}
Thus, we see that $\tilde{E}(x,s)=(\ML E)(x,s)$ is absolutely integrable on vertical lines $\Re s= a$, $a>0$, so $E(x,t) = \frac{1}{2\pi\I}\int_{a-\I\infty}^{a+\I\infty}\tilde{E}(x,s)\e^{ts}\dif s$ is a continuous function of $x$ and $t$.
\end{proof}

Actually one can show that 
\begin{align*}
	\WF(S) = {}	&\{(0,0;\xi,\eta): (\xi, \eta)\neq(0,0)\} \cup\{(0,t; \xi,0): t>0, \xi\neq0\} \\
		{}\cup{} 	&\{(x, t; \xi, \eta): t>0, (\xi,\eta)\neq(0,0),\ \abs{x}=t/\sqrt{\tau}, (x,t)\cdot(\xi,\eta)=0\}.
\end{align*}
For $\abs{x} < t/\sqrt{\tau}$, one shifts the contour in the inverse Laplace transform to the left to obtain 
\[
	S(x,t) = \frac{1}{2} + \frac{1}{4\pi\I}\int_{\Gamma}\frac{l_{1}(s)}{s}\exp\bigl(-\abs{x}sl_{1}(s) + ts\bigr)\dif s, 
\]
where $\Gamma$ is a (finite) contour encircling the branch cut $[-1/\tau, -1]$ in the counterclockwise direction. This shows that $S$ is analytic on the set $0<\abs{x}<t/\sqrt{\tau}$. By differentiating the above formula with respect to $x$ and using the residue theorem, one gets that 
\[
	\dpd{S}{x}(0^{+}, t) = -\frac{1-\tau}{2}\e^{-t}, \quad \dpd{S}{x}(0^{-},t) = \frac{1-\tau}{2}\e^{-t}, \quad t>0.
\]
As in the proof of Theorem \ref{th: Gevrey regularity}, $S$ is analytic in the $(0,t)$-direction at points $(0,t_{0})$, $t_{0}>0$. For $x\neq0$, consider $S^{(-1)}(x,t) = \int_{0}^{t}S(x,t_{1})\dif t_{1}$. Changing variables as in the proof of Theorem \ref{th: Gevrey regularity} and differentiating with respect to $u=\sqrt{\tau}x+ t$, we get
\[
	\dpd[n]{(S^{(-1)})^{\natural}}{u}(u,v) = \frac{1}{4\pi\I}\int_{a-\I\infty}^{a+\I\infty}\frac{l_{1}(s)}{s^{2}}\Bigl\{\frac{s}{2}\Bigl(1-\frac{l_{1}(s)}{\sqrt{\tau}}\Bigr)\Bigr\}^{n}\exp\Bigl\{\frac{us}{2}\Bigl(1-\frac{l_{1}(s)}{\sqrt{\tau}}\Bigr) - \frac{vs}{2}\Bigl(1+\frac{l_{1}(s)}{\sqrt{\tau}}\Bigr)\Bigr\}\dif s.
\]
In view of \eqref{eq: approx l_{1}}, this integral converges absolutely for every $n$, and is bounded even by $D^{n+1}$ for some $D>0$, showing the $(S^{(-1)})^{\natural}$ is real analytic in the $(u,0)$-direction at points $(u_{0},0)$, $u_{0}>0$. We omit further details.

\mbox{}

We can again investigate the response in this model to a forced oscillation at the origin, starting at $t=0$. As before, we set $f(x,t)=\delta(x)H(t)\cos(\omega t)$, with $\omega>0$, and $u_{0}=v_{0}=0$. Then the solution $u$ of the initial value problem \eqref{eq: Cauchy problem Zener} has Laplace transform
\[
	\tilde{u}(x,s) = \frac{l_{1}(s)}{2}\exp\bigl(-\abs{x}sl_{1}(s)\bigr)\frac{1}{s^{2}+\omega^{2}}.
\]
Notice that this Laplace transform is integrable on vertical lines $\Re s=a$, so $u$ is continuous (although it is not of class $C^{1}$). Again the solution has support inside the cone $t\ge \sqrt{\tau}\abs{x}$. For $x$ and $t$ with $t > \sqrt{\tau}\abs{x}$, we can move the contour in the Inverse Laplace transform to the left to get
\[
	u(x,t) = H(t/\sqrt{\tau} - \abs{x}) \bigl(u_{\mathrm{ss}}(x,t) + u_{\mathrm{ts}}(x,t)\bigr),
\]
with 
\begin{align*}
	u_{\mathrm{ss}}(x,t)	&= \frac{l_{1}(\I\omega)}{4\I\omega}\exp\bigl(-\abs{x}\I\omega l_{1}(\I\omega) + \I\omega t\bigr) 
							- \frac{l_{1}(-\I\omega)}{4\I\omega}\exp\bigl(\,\abs{x}\I\omega l_{1}(-\I\omega) - \I\omega t\bigr) \\
						&=\frac{\rho_{1}(\omega)}{2\omega}\e^{-b_{1}(\omega)\omega\,\abs{x}}\sin\bigl(\omega t - a_{1}(\omega)\omega \abs{x} - \phi_{1}(\omega)\bigr); \\	
	u_{\mathrm{ts}}(x,t)		&= \frac{1}{4\pi\I}\int_{\Gamma}l_{1}(s)\exp\bigl(-\abs{x}sl_{1}(s) + ts\bigr)\frac{1}{s^{2}+\omega^{2}}\dif s;
\end{align*}
where again $\Gamma$ is a (finite) closed contour encircling the branch cut $[-1/\tau, -1]$ in the counterclockwise direction, and $l_{1}(\I\omega) = {a_{1}(\omega) - \I b_{1}(\omega)} = \rho_{1}(\omega)\e^{-\I\phi_{1}(\omega)}$, as before. The transient state $u_{\mathrm{ts}}$ converges to $0$ as $t\to\infty$, locally uniformly in $x$. For fixed $\omega$, the steady state $u_{\mathrm{ss}}$ is formally identical to the steady states in the fractional Zener model \eqref{eq: steady state}. We have the complex dispersion relation $k_{1}(\omega) = \omega l_{1}(\I\omega)$, and phase velocity $V_{1}(\omega) = 1/a_{1}(\omega)$. However, there is a qualitative difference between the two models in terms of the dependency of the dissipation on the frequency $\omega$. In the SLS model (corresponding to $\alpha=1$), the attenuation coefficient $d_{1}(\omega) =  b_{1}(\omega)\omega$ has asymptotic behavior 
\[
	d_{1}(\omega) \sim \frac{\sqrt{\tau}}{2}\Bigl(\frac{1}{\tau}-1\Bigr), \quad \text{as } \omega\to\infty,
\]
contrasting \eqref{eq: asymp behavior delta(omega)}, which shows that the attenuation coefficient in the fractional Zener model ($0<\alpha<1$) grows to $\infty$ as $\omega\to\infty$. In the SLS model, two pseudo-monochromatic waves with different frequencies have roughly the same amount of spatial dampening, while in the fractional Zener model, the wave with the higher frequency will experience more dampening than the wave with the lower frequency.

%%%%%%%%%%%%%%%%%%%%%%%%%%%%%%%%%%%%%%%%%%%%%%%%%%%%%%%%%%%%%%%%%%%%%%%%%%%%%%
%% DISTRIBUTED-ORDER FRACTIONAL WAVE EQUATIONS %%%%%%%%%%%%%%%%%%%%%%%%%%%%%%%%%%%%%%%%%%%%%%
%%%%%%%%%%%%%%%%%%%%%%%%%%%%%%%%%%%%%%%%%%%%%%%%%%%%%%%%%%%%%%%%%%%%%%%%%%%%%%
%\section{Distributed-order fractional wave equations}
%
%\subsection{Existence and uniqueness of solutions}
%
%\begin{theorem}
%Let $u_{0}, v_{0} \in \mathcal{E}'(\R)$. Suppose $\dme$ and $\dms$ be positive Radon measures on $[0,1]$, satisfying condition \eqref{T}. Then the equation \eqref{...} has a unique solution $u\in ...$, given by
%\[
%	u(x,t) = 
%\] 
%\end{theorem}
%
%\subsection{Micro-local analysis}

\appendix
%%%%%%%%%%%%%%%%%%%%%%%%%%%%%%%%%%%%%%%%%%%%%%%%%%%%%%%%%%%%%%%%%%%%%%%%%%%%%
% APPENDIX %%%%%%%%%%%%%%%%%%%%%%%%%%%%%%%%%%%%%%%%%%%%%%%%%%%%%%%%%%%%%%%%%%%%%%
%%%%%%%%%%%%%%%%%%%%%%%%%%%%%%%%%%%%%%%%%%%%%%%%%%%%%%%%%%%%%%%%%%%%%%%%%%%%%
\section{Proof of Lemma \ref{lem: technical lemma}}
\label{sec: technical calculation}
In this appendix we provide a proof Lemma \ref{lem: technical lemma}. We use the same notations as in the proof of Theorem \ref{th: Gevrey regularity}. Set
\begin{align*}
	h_{m}(w)	&\coloneqq f(w) + g_{m}(w), \quad \text{where}\\
	f(w)		&= \kappa w - \frac{1}{1-\alpha}w^{1-\alpha} + \log w + \log 2, \\
	g_{m}(w) 	&= \frac{\kappa w}{2}E_{1}(\mu m^{\frac{1}{1-\alpha}} w) - \frac{u\mu m^{\frac{\alpha}{1-\alpha}}w}{2}E_{2}(\mu m^{\frac{1}{1-\alpha}} w) 
				+ \log\Bigl(1+\frac{E_{1}(\mu m^{\frac{1}{1-\alpha}} w)}{2}\Bigr).
\end{align*}
Here, $\mu$ is a fixed constant, $E_{1}$ and $E_{2}$ are ``remainder functions'' given by \eqref{eq: remainders}, and $\kappa$ is a number in a fixed range, namely
\begin{equation}
\label{eq: range kappa}
	\kappa \in I \coloneqq \bigl[\bigl(1000/\sin(\alpha\pi)\bigr)^{\frac{1}{1-\alpha}}, \bigl(2000/\sin(\alpha\pi)\bigr)^{\frac{1}{1-\alpha}}\bigr].
\end{equation}
We need to show that we can choose $\kappa=\kappa_{m}$ in the interval $I$ in such a way that
\begin{equation}
\label{eq: bound im integral}	
	\Im \int_{0}^{+\I\infty}l_{\alpha}(\mu m^{\frac{1}{1-\alpha}}w)\e^{mh_{m}(w)}\frac{\dif w}{w} \gtrsim \frac{c^{m}}{\sqrt{m}}. 
\end{equation}
for some $c>0$ independent from $m$.

First we show necessary estimates uniformly for $\kappa\in I$, and later we show the existence of $\kappa_{m}$ so that \eqref{eq: bound im integral} holds.
In order to prove \eqref{eq: bound im integral}, we will use the saddle point method. We view the phase function $h_{m}$ as a small perturbation of $f$. Indeed, by the bounds \eqref{eq: bounds E1 and E2} we get 
\begin{equation}
\label{eq: asymptotic g_{m}}
	g_{m}(w)\lesssim m^{-\frac{\alpha}{1-\alpha}}\abs{w}^{1-\alpha}, \quad m\to\infty,
\end{equation}
provided that $\abs{w}$ is greater than some fixed $\eps>0$. In particular, $g_{m}$, as well as its derivatives, converges locally uniformly to $0$ on the set $\abs{w}>\eps$. This convergence is also uniform with respect to $\kappa\in I$. We want the quantity $g_{m}$ and its derivatives to be very small near the saddle point, so we also bound $\abs{w}$ from above, and consider the fixed range $\eps \le \abs{w} \le W$. We will perform the saddle point analysis in this range, and show that the parts of the integral with $\abs{w}\le \eps$ and $\abs{w}\ge W$ are negligible with respect to the contribution from the saddle point. We now set
\begin{equation}
\label{eq: def eps and W}
	\eps = \frac{1}{1+1/\sqrt{\tau}}\biggl(\frac{\sin(\alpha\pi)}{2000}\biggr)^{\frac{4}{1-\alpha}}, \quad W = \frac{8}{1-\alpha}.
\end{equation}
%Although the final contour will lie in the quadrant $\arg w \in [\pi/2, \pi]$, 
It will be useful to consider $f$ and $g_{m}$ as holomorphic functions on the set $\Omega = \{ w: \eps\le\abs{w}\le W, 0\le \arg w\le 3\pi/2\}$ by analytic continuation\footnote{Note that this continuation is different from the one to $\arg s \in [-\pi,-\pi/2]$, which appeared in (the derivation of) \eqref{eq: S in cone}. They are situated on different sheets of the Riemann surface of the logarithm.}. If $m$ is sufficiently large, then the zero and the pole of the function $(1+\tau(\mu m^{\frac{1}{1-\alpha}} w)^{\alpha})/(1+(\mu m^{\frac{1}{1-\alpha}} w)^{\alpha})$ have modulus strictly smaller than $\eps$, so that $g_{m}$ is well defined and holomorphic in $\Omega$. We let $m$ be so large that 
\begin{equation}
\label{eq: bound g}
	\abs[1]{g^{(i)}_{m}(w)} < \frac{\sin^{2}(\alpha\pi)}{1000^{2}}, \quad \text{for } w\in\Omega, \quad \kappa\in I, \quad i=0,1,2,3.
\end{equation}

Next we focus our attention on $f$. We have 
\begin{align*}
	f'(w)	&= \kappa - w^{-\alpha} + \frac{1}{w}. 
\end{align*}
To solve the saddle point equation $f'(w)=0$, it is convenient to solve for $z\coloneqq 1/w$. The equation $\kappa - z^{\alpha} + z = 0$ will have a solution $z_{0}$ near $\kappa\e^{-\I\pi}$. Indeed, noting that $\abs{z^{\alpha}} < \abs{\kappa+z}$ on $\partial B(\kappa\e^{-\I\pi}, \kappa/2)$, it follows from Rouch\'e's theorem that $\kappa - z^{\alpha} + z$ has a unique zero in the disc $B(\kappa\e^{-\I\pi}, \kappa/2)$. Next, we will deduce a precise estimate for $z_{0}$. We want to keep all the constants in the error terms explicit. For this, we will use $\zeta$ to denote some complex number with $\abs{\zeta}\le 1$. At each next occurrence, this number might have a different value than the previous occurrence, but we will use the same notation $\zeta$ each time. We have
\begin{align*}
	z_{0}	&= \e^{-\I\pi}\kappa\frac{1}{1-z_{0}^{\alpha-1}} = \e^{-\I\pi}\kappa\bigl(1+z_{0}^{\alpha-1} + 8\zeta\abs{z_{0}}^{2\alpha-2}\bigr) 
					= \e^{-\I\pi}\kappa\Bigl(1+z_{0}^{\alpha-1} + \frac{16\sin^{2}(\alpha\pi)}{1000^{2}}\zeta\Bigr) \\
		&= \e^{-\I\pi}\kappa\Bigl(1 + \e^{\I\pi(1-\alpha)}\kappa^{\alpha-1}\Bigl(1+\frac{6\sin(\alpha\pi)}{1000}\zeta\Bigr) + \frac{16\sin^{2}(\alpha\pi)}{1000^{2}}\zeta\Bigr) \\	
		&= \e^{-\I\pi}\kappa\Bigl(1 + \e^{\I\pi(1-\alpha)}\kappa^{\alpha-1} + \frac{22\sin^{2}(\alpha\pi)}{1000^{2}}\zeta\Bigr).		
\end{align*}
Here we used Taylor's theorem with explicit error terms, the a priori estimate $\abs{z_{0}} > \kappa/2$, and the bound $\kappa^{\frac{1}{1-\alpha}} \ge 1000/\sin(\alpha\pi)$. For $w_{0}=1/z_{0}$ we get
\begin{align}
	w_{0} 		&= \frac{\e^{\I\pi}}{\kappa}\Bigl(1-\e^{\I\pi(1-\alpha)}\kappa^{\alpha-1} + \frac{54\sin^{2}(\alpha\pi)}{1000^{2}}\zeta\Bigr), \label{eq: approx w_{0}}\\
	\arg(w_{0})	&= \pi - \frac{2\sin(\alpha\pi)}{998}\xi, \quad \text{for some } 0<\xi<1. \nonumber
\end{align}
Let us denote by $w_{m}$ the saddle point of $h_{m}= f+g_{m}$. By Hurwitz's theorem, we may assume that
\begin{equation}
\label{eq: approx w_{m}}
w_{m} = w_{0}\Bigl(1+\frac{\sin^{2}(\alpha\pi)}{1000^{2}}\zeta\Bigr) = \frac{\e^{\I\pi}}{\kappa}\Bigl(1-\e^{\I\pi(1-\alpha)}\kappa^{\alpha-1} + \frac{56\sin^{2}(\alpha\pi)}{1000^{2}}\zeta\Bigr)
\end{equation}
for sufficiently large $m$.
%$= (\e^{\I\pi}/\kappa)(1+(4\alpha/1000)\zeta)$ say.
We will let the contour pass through $w_{m}$ via the steepest path. 

\begin{lemma}
\label{lem: steepest path}
There exists some $\delta>0$, and a contour $\Gamma$, the path of steepest descent, which connects two (nearly) opposing points $c$ and $d$ on the circle $\abs{w-w_{m}}=\delta$. This path passes through $w_{m}$, $\Im h_{m}$ is constant along it, while $\Re h_{m}$ reaches its maximum at $w_{m}$. The tangent vector along $\Gamma$ has its argument in the range $(3\pi/4, 5\pi/4)$.
\end{lemma}
The proof will show that we may take $\delta= (27/680)\kappa^{-1}$.
\begin{proof}
The idea is to approximate $h_{m}(w)-h_{m}(w_{m})$ by the quadratic function ${(h_{m}''(w_{m})/2)(w-w_{m})^{2}}$. By Taylor's theorem, we have on a small neighborhood of $w_{m}$
\[
	h_{m}(w) - h_{m}(w_{m}) = \frac{h_{m}''(w_{m})}{2}(w-w_{m})^{2}(1+\eta_{m}(w)),
\]
where $\eta_{m}(w)$ is a holomorphic function satisfying 
\[
	\abs{\eta_{m}(w)} \le \biggl(\frac{2}{3!\abs{h_{m}''(w_{m})}}\max_{z\in[w_{m}, w]}\abs{h_{m}'''(z)}\biggr)\abs{w-w_{m}}.
\]
The derivatives of $h_{m}$ can be approximated by those of $f$. We have
\[
	f''(w) =  -\frac{1}{w^{2}}\bigl(1-\alpha w^{1-\alpha}\bigr), \quad f'''(w) = \frac{2}{w^{3}}\bigl(1-\alpha(\alpha+1)w^{1-\alpha}\bigr),
\]
so that by \eqref{eq: bound g} and \eqref{eq: approx w_{m}},  
\[
	\abs{h_{m}''(w_{m})} 	\ge \abs{f''(w_{m})} - \frac{\sin^{2}(\alpha\pi)}{1000^{2}} \ge \kappa^{2}\Bigl(\frac{998}{1000}\Bigr)^{2}\Bigl(1-\frac{2}{1000}\Bigr) - \frac{\sin^{2}(\alpha\pi)}{1000^{2}} \ge \frac{9\kappa^{2}}{10}.
\]
Also for $\abs{w-\e^{\I\pi}/\kappa} \le 1/(2\kappa)$, 
\[
	\abs{h_{m}'''(w)} \le \frac{2}{1/(2\kappa)^{3}}\bigl(1+2(3/(2\kappa))^{1-\alpha}\bigr) + \frac{\sin^{2}(\alpha\pi)}{1000^{2}} \le 17\kappa^{3}.
\]
Hence we get 
\begin{equation}
\label{eq: bound eta}
\abs{\eta_{m}(w)} \le (170/27)\kappa\abs{w-w_{m}}, \quad \text{for } \abs[1]{w-\e^{\I\pi}/\kappa} \le (1/2)\kappa^{-1}.
\end{equation}
We now set $\delta = (27/680)\kappa^{-1}$, then for $\abs{w-w_{m}} \le \delta$, 
\begin{align}
	h_{m}(w) - h_{m}(w_{m}) 	&= \frac{h_{m}''(w_{m})}{2}(w-w_{m})^{2}(1+\eta_{m}(w)), \quad \abs{\eta_{m}(w)} \le \frac{1}{4}, \label{eq: h near saddle point} \\
						&\eqqcolon \bigl(\psi_{m}(w)\bigr)^{2}. \nonumber
\end{align}
Here\footnote{The minus sign here is introduced for convenience. With this minus sign, the steepest path defined later on will have the desired orientation.}, $\psi_{m}(w) = -\sqrt{h_{m}''(w_{m})/2}(w-w_{m})\sqrt{1+\eta_{m}(w)}$, where $\sqrt{\phantom{a}}$ denotes the principal branch of the square root. We claim that this is a holomorphic bijection from the closed disk $\overline{B}(w_{m},\delta)$ onto some compact neighborhood $F$ of zero. This will follow if we show that its derivative does not vary too much. We have 
\[
	\psi_{m}'(w) = -\sqrt{\frac{h_{m}''(w_{m})}{2}}\biggl(\sqrt{1+\eta_{m}(w)} + (w-w_{m})\frac{\eta_{m}'(w)}{2\sqrt{1+\eta_{m}(w)}}\biggr).
\]
Estimating $\eta_{m}'$ by Cauchy's formula, if $\abs{w-w_{m}}\le\delta$, then
\[
	\abs{\eta_{m}'(w)} \le \frac{1}{2\pi}\abs{\int_{\partial B(w_{m}, 12\delta)}\frac{\eta_{m}(z)}{(z-w)^{2}}\dif z} \le \frac{1}{2\pi}\cdot\frac{1}{(11\delta)^{2}}\cdot 3 \cdot (24\delta\pi) = \frac{36}{121\delta}.
\] 
Here we used that $\abs{z-w}\ge11\delta$, and that $\abs{\eta_{m}(z)} \le 3$ on the circle $\abs{z-w_{m}} = 12\delta$, by \eqref{eq: bound eta}. Hence, if $\abs{w-w_{m}}\le\delta$, then
\[
	\psi_{m}'(w) = -\sqrt{\frac{h_{m}''(w_{m})}{2}} \Bigl(1+\frac{\zeta_{w}}{2}\Bigr), \quad \text{for some $\zeta_{w}$ with $\abs{\zeta_{w}}\le1$}.
\]
In particular, if $w, w' \in \overline{B}(w_{m}, \delta)$, $w\neq w'$, then $\psi_{m}(w')-\psi_{m}(w) = \int_{w}^{w'}\psi_{m}'(z)\dif z \neq 0$, showing that $\psi_{m}$ is injective.

We now set $\Gamma \coloneqq \psi_{m}^{-1}(L)$, where $L = [\I a, \I b]$ is the maximal line segment along the imaginary axis, which contains the origin and lies completely within $F$. This line segment connects the boundary points $\I a, \I b\in \partial F$ via the imaginary axis, passing through the origin, and $\Gamma$ is a path which connects the points $c\coloneqq\psi_{m}^{-1}(\I a)$ and $d\coloneqq \psi_{m}^{-1}(\I b)$ on the circle $\partial B(w_{m}, \delta)$ via a path passing through the saddle point $w_{m}$. For points $w\in\Gamma$, clearly $\Im h_{m}(w)=\Im h_{m}(w_{m})$, and $\Re h_{m}(w)\le\Re h_{m}(w_{m})$, with equality only when $w=w_{m}$.

Let us now locate the points $c$ and $d$ on the circle a bit more, using \eqref{eq: h near saddle point}. Writing $h_{m}''(w_{m}) = \abs{h_{m}''(w_{m})}\e^{\I\vphi_{m}}$, similar calculations as before show that 
\[
	\vphi_{m} = \pi - \frac{8\xi}{992}, \quad \text{for some $\xi$ with $0 < \xi < 1$.}
\]
Setting $c=w_{m}+\delta\e^{\I\theta_{c}}$, $d=w_{m}+\delta\e^{\I\theta_{d}}$, we get the following relations for $\theta = \theta_{c}, \theta_{d}$, by taking real and imaginary part in \eqref{eq: h near saddle point}:
\begin{align*}
	0	&> \frac{\abs{h_{m}''(w_{m})}}{2}\delta^{2}\bigl(	\cos(\vphi_{m}+2\theta)(1+\Re\eta_{m}(\delta\e^{\I\theta})) - \sin(\vphi_{m}+2\theta)\Im\eta_{m}(\delta\e^{\I\theta})\bigr), \\
	0	&= \frac{\abs{h_{m}''(w_{m})}}{2}\delta^{2}\bigl( \cos(\vphi_{m}+2\theta)\Im\eta_{m}(\delta\e^{\I\theta}) + \sin(\vphi_{m}+2\theta)(1+\Re\eta_{m}(\delta\e^{\I\theta}))\bigr).
\end{align*}
The first inequality implies that $\theta \in (-3\pi/8, 3\pi/8) \cup (5\pi/8, 11\pi/8)$ (if not, then using the estimate of $\vphi_{m}$ and the bound on $\eta_{m}$, one shows that the right-hand side would be positive). Using this initial localization, we can narrow the range down using the second equality, to $\theta \in (-\pi/8, \pi/8) \cup (7\pi/8, 9\pi/8)$ say. Actually, $\theta_{c}\in (-\pi/8, \pi/8)$, while $\theta_{d}\in(7\pi/8, 9\pi/8)$. This will follow from the following estimate for the argument of the tangent vector along $\Gamma$.

We have the following parametrization of $\Gamma$: $\gamma: [a, b] \to \Gamma: y \mapsto \psi_{m}^{-1}(\I y)$. We have
\[
	\gamma'(y) = \frac{\I}{\psi_{m}'(\psi_{m}^{-1}(\I y))} = -\I \sqrt{\frac{2}{h_{m}''(w_{m})}}\Bigl(1+\frac{\zeta_{y}}{2}\Bigr),
\]
for some $\zeta_{y}$ with $\abs{\zeta_{y}}\le1$. Using the bounds on $\vphi_{m}$, we see that this tangent vector has its argument in the range $(5\pi/6 - 4/992, 7\pi/6 + 4/992)$.
\end{proof}

We will now use the obtained information to estimate $\int_{\Gamma}l_{\alpha}(\mu m^{\frac{1}{1-\alpha}}w)w^{-1}\e^{mh_{m}(w)}\dif w$, with $\Gamma$ as in the above lemma. First we reparametrize $\Gamma$ with arc length:
\[
	\tilde{\gamma}: [\tilde{a}, \tilde{b}] \to \Gamma: u \mapsto \tilde{\gamma}(u), \quad \tilde{\gamma}(0)=w_{m}, \quad \abs{\tilde{\gamma}'(u)} = 1.
\]	
From the lemma, we have $\abs{\arg\tilde{\gamma}'(u)-\pi} < \pi/4$. If $m$ is sufficiently large, then $\abs[1]{\arg l_{\alpha}(\mu m^{\frac{1}{1-\alpha}}w)} \le 1/1000$ and $\abs{\arg w^{-1} + \pi} \le \arctan \frac{2/1000+27/680}{998/1000-27/680} \le 1/20$, for $w$ on the steepest path $\Gamma$. We have 
\begin{align*}
	\int_{\Gamma}l_{\alpha}(\mu m^{\frac{1}{1-\alpha}}w)\e^{mh_{m}(w)}\frac{\dif w}{w} 	
			&= \e^{mh_{m}(w_{m})}\int_{\tilde{a}}^{\tilde{b}}l_{\alpha}(\mu m^{\frac{1}{1-\alpha}}\tilde{\gamma}(u))
				\exp\bigl\{m\bigl(h_{m}(\tilde{\gamma}(u))-h_{m}(w_{m})\bigr)\bigr\}\frac{\tilde{\gamma}'(u)}{\tilde{\gamma}(u)}\dif u \\
			&\eqqcolon \e^{mh_{m}(w_{m})}R\e^{\I \phi}.	
\end{align*}
Here, $\abs{\phi} \le \pi/4 + 1/1000 + 1/20$, and $R \ge \cos(\pi/4 + 1/1000 + 1/20)\int_{\tilde{a}}^{\tilde{b}}\abs{\dotso} \dif u$. Using \eqref{eq: h near saddle point} and the estimates $l_{\alpha} \gtrsim 1$, $1/\tilde{\gamma}(u) \gtrsim \kappa$, and $\abs{h''_{m}(w_{m})} \simeq \kappa^{2}$, we get
\begin{align*}
	R 	&\gtrsim \kappa \int_{\tilde{a}}^{\tilde{b}}\exp\biggl\{m\biggl(\frac{h_{m}''(w_{m})}{2}(\tilde{\gamma}(u)-w_{m})^{2}(1+\eta_{m}(\tilde{\gamma}(u)))\biggr)\biggr\}\dif u \\
		&\gtrsim \kappa\int_{\tilde{a}}^{\tilde{b}}\exp\biggl\{-m\frac{5\abs{h_{m}''(w_{m})}}{8}u^{2}\biggr\}\dif u \gtrsim \kappa\cdot \frac{1}{\sqrt{m\abs{h_{m}''(w_{m})}}} \gtrsim \frac{1}{\sqrt{m}}, \quad \text{as $m\to\infty$}.
\end{align*}
Here we used that the exponent is real and non-positive along $\Gamma$, so that 
\begin{align*}
	\frac{h_{m}''(w_{m})}{2}(\tilde{\gamma}(u)-w_{m})^{2}(1+\eta_{m}(\tilde{\gamma}(u))) 	&= -\frac{\abs{h_{m}''(w_{m})}}{2}\abs{\tilde{\gamma}(u)-w_{m}}^{2}\abs{1+\eta_{m}(\tilde{\gamma}(u))} \\
																		&\ge -\frac{\abs{h_{m}''(w_{m})}}{2}\cdot u^{2} \cdot \frac{5}{4}.	
\end{align*}
The last inequality follows from the fact that $u$ is arc length: $\abs{\tilde{\gamma}(u)-w_{m}} = \abs{\tilde{\gamma}(u)-\tilde{\gamma}(0)} \le \abs{u}$. For the same reason $\abs[0]{\tilde{b}-\tilde{a}} = \length(\Gamma) \gtrsim 1/\kappa$, so that the integral above can be transformed to the integral of $\e^{-t^{2}}$ over an interval containing $0$ of length $\gtrsim 1$. 
We can conclude that 
\begin{align}
	\abs{\int_{\Gamma}l_{\alpha}(\mu m^{\frac{1}{1-\alpha}}w)\e^{mh_{m}(w)}\frac{\dif w}{w}} 	
		& \gtrsim \frac{\e^{m\Re h_{m}(w_{m})}}{\sqrt{m}} \label{eq: contribution saddle point}\\
		& \gtrsim \frac{c^{m}}{\sqrt{m}},  \quad c = \biggl(\frac{\sin(\alpha\pi)}{2000}\biggr)^{\frac{2}{1-\alpha}}. \label{eq: def c} 
\end{align}
The value for $c$ arises from the following (rough) lower bound for $\Re h_{m}(w_{m})$: using $1/(2\kappa) \le \abs{w_{m}} \le 2/\kappa$ and $\kappa\in I$ we get
\begin{align*}
	\Re h_{m}(w_{m}) 	&= \Re\Bigl(\kappa w_{m} -\frac{1}{1-\alpha}w_{m}^{1-\alpha} + \log w_{m} + \log 2 + g_{m}(w_{m})\Bigr) \\
					&\ge -2 - \frac{2\sin(\alpha\pi)}{1000(1-\alpha)} - \frac{1}{1-\alpha}\log\frac{2000}{\sin(\alpha\pi)} - \frac{1}{1000^{2}}\\
					&\ge -\frac{2}{1-\alpha}\log\frac{2000}{\sin(\alpha\pi)}.
\end{align*}

To control the phase of $\int_{\Gamma}$, we need a precise estimate of $\Im h_{m}(w_{m})$. Using \eqref{eq: approx w_{m}}, we get (now using $\xi$ for a \emph{real} number satisfying $\abs{\xi}\le1$, with a possibly different value at each occurrence)
\begin{align*}
	\Im w_{m} 		&= \frac{1}{\kappa}\Bigl(\sin(\alpha\pi)\kappa^{\alpha-1} + \frac{56\sin^{2}(\alpha\pi)}{1000^{2}}\xi\Bigr); \\
	\Im w_{m}^{1-\alpha}	&= \sin(\alpha\pi)\kappa^{\alpha-1}\Bigl(1+\frac{4\sin(\alpha\pi)}{1000}\xi\Bigr) = \sin(\alpha\pi)\kappa^{\alpha-1}+\frac{4\sin^{2}(\alpha\pi)}{1000^{2}}\xi;\\
	\arg w_{m}		&= \pi - \arctan\biggl(\frac{\sin(\alpha\pi)\kappa^{\alpha-1}+ \frac{56\sin^{2}(\alpha\pi)}{1000^{2}}\xi}{1+\frac{2\sin(\alpha\pi)}{1000}\xi}\biggr) 
					= \pi - \sin(\alpha\pi)\kappa^{\alpha-1} + \frac{63\sin^{2}(\alpha\pi)}{1000^{2}}\xi.				
\end{align*}
This implies that 
\begin{align*}
	\Im h_{m}(w_{m}) 	&= \kappa\Im w_{m} - \frac{1}{1-\alpha}\Im w_{m}^{1-\alpha} + \arg w_{m} + \Im g_{m}(w_{m}) \\
					&= \pi -\frac{\sin(\alpha\pi)}{1-\alpha}\Bigl(\kappa^{\alpha-1} + \frac{124\sin(\alpha\pi)}{1000^{2}}\xi\Bigr),
\end{align*}
where we also used \eqref{eq: bound g} to bound $\Im g_{m}(w_{m})$. Now that we have such a precise estimate for $\Im h_{m}(w_{m})$, we will demonstrate how to choose $\kappa\in I$. We have 
\begin{align*}
	&\bigl(\Im h_{m}(w_{m})\bigr)\big\rvert_{\kappa^{1-\alpha} = \frac{1000}{\sin(\alpha\pi)}} = \pi -\frac{\sin(\alpha\pi)}{1-\alpha}\Bigl(\frac{\sin(\alpha\pi)}{1000} + \frac{124\sin(\alpha\pi)}{1000^{2}}\xi\Bigr) \\
<{}	&\bigl(\Im h_{m}(w_{m})\bigr)\big\rvert_{\kappa^{1-\alpha} = \frac{2000}{\sin(\alpha\pi)}} = \pi -\frac{\sin(\alpha\pi)}{1-\alpha}\Bigl(\frac{\sin(\alpha\pi)}{2000} + \frac{124\sin(\alpha\pi)}{1000^{2}}\xi'\Bigr).
\end{align*}
We now note that the value of $w_{m}$ depends continuously on $\kappa$, and so also $h_{m}(w_{m})$ depends continuously on $\kappa$. Hence, for each sufficiently large $m$, we may choose $\kappa=\kappa_{m} \in I$ in such a way that $m\bigl(\Im h_{m}(w_{m})\bigr) \in \pi/2 + 2\pi\Z$. This guarantees that 
\[
	\Im \int_{\Gamma}l_{\alpha}(\mu m^{\frac{1}{1-\alpha}}w)\e^{mh_{m}(w)}\frac{\dif w}{w} = \Im \bigl(\e^{mh_{m}(w_{m})}R\e^{\I \phi}\bigr) \gtrsim \frac{c^{m}}{\sqrt{m}}.
\]

Finally, we have to deform the complete contour $[0, +\I\infty)$ to a contour containing $\Gamma$, and show that the contribution from the other pieces is negligible. We do this in several steps. First, we set $\Upsilon_{1} = [0, \I\eps]$. For points $w$ with $\abs{w} < \eps$, the asymptotic estimates \eqref{eq: bounds E1 and E2} on the remainder functions $E_{1}$ and $E_{2}$ cannot be used. Writing the integrand in its original form, that is before introducing the functions $E_{1}$, $E_{2}$, $f$, and $g$, we get 
\begin{align*}
	\int_{\Upsilon_{1}}l_{\alpha}(\mu m^{\frac{1}{1-\alpha}}w)\e^{mh_{m}(w)}\frac{\dif w}{w} 	
	&= \int_{\Upsilon_{1}}l_{\alpha}(\mu m^{\frac{1}{1-\alpha}}w)w^{m}\biggl(\frac{l_{\alpha}(\mu m^{\frac{1}{1-\alpha}}w)}{\sqrt{\tau}} + 1\biggr)^{m}\\
	&\exp\biggl\{-\frac{u\mu m^{\frac{1}{1-\alpha}}w}{2}\biggl(\frac{l_{\alpha}(\mu m^{\frac{1}{1-\alpha}}w)}{\sqrt{\tau}} - 1\biggr) + \frac{m\kappa w}{2}\biggl(\frac{l_{\alpha}(\mu m^{\frac{1}{1-\alpha}}w)}{\sqrt{\tau}} + 1\biggr)\biggr\}\frac{\dif w}{w} \\
	&\lesssim \eps^{m}\Bigl(\frac{1}{\sqrt{\tau}}+1\Bigr)^{m} = c^{2m}.	
\end{align*}
Here we used that $\abs{l_{\alpha}(s)} \le 1$ and $\Im l_{\alpha}(s) \le 0$ for $s\in \I\R_{+}$, and the definitions of $\eps$ and $c$, \eqref{eq: def eps and W} and \eqref{eq: def c}. Since $c<1$, this is of strictly lower order than the contribution from the integral over $\Gamma$.

Next we set $\Upsilon_{2}\coloneqq\{\eps \e^{\I\vphi}: \pi/2\le \vphi \le \theta_{m}\}$, were $\theta_{m}=\arg w_{m}$. We have
\begin{align*}
	\Re h_{m}(\eps\e^{\I\vphi})		&= \kappa\eps\cos\vphi - \frac{1}{1-\alpha}\eps^{1-\alpha}\cos((1-\alpha)\vphi) + \log\eps + \log 2 + \Re g(\eps\e^{\I\vphi}) \\
							&\le \frac{1}{1-\alpha} + \frac{1}{1000^{2}} -\frac{4}{1-\alpha}\log\frac{2000}{\sin(\alpha\pi)} \le \frac{3}{2}\log c. 
\end{align*}
Hence,
\[
	\int_{\Upsilon_{2}}l_{\alpha}(\mu m^{\frac{1}{1-\alpha}}w)\e^{mh_{m}(w)}\frac{\dif w}{w} \lesssim c^{3m/2},
\]
which is negligible.

Next we set $\Upsilon_{3} \coloneqq \{r\e^{\I\theta_{m}}: \eps\le r \le r_{1}\}$, where $r_{1}$ is such that this line segment connects $\eps\e^{\I\theta_{m}}$ to the circle $\abs{w-w_{m}} = \delta$, so $r_{1} = \abs{w_{m}} - \delta$. Note that $r_{1}\approx \frac{653}{680}\kappa^{-1} \ge \frac{1}{2}\bigl(\sin(\alpha\pi)/2000\bigr)^{\frac{1}{1-\alpha}} > \eps$. Consider the function $r \mapsto \Re f(r\e^{\I\theta_{m}})$. This function is non-decreasing for $r\in[\eps, r_{1}]$. Indeed, using that $r\le r_{1} \le \frac{682-27}{680}\kappa^{-1}$ we get
\begin{align*}
	\dpd{}{r}\Re f(r\e^{\I\theta_{m}}) 	&= \dpd{}{r}\biggl(\kappa r \cos(\theta_{m}) - \frac{1}{1-\alpha}r^{1-\alpha}\cos((1-\alpha)\theta_{m}) + \log r + \log 2\biggr)\\
								&= \kappa\cos(\theta_{m}) - \cos((1-\alpha)\theta_{m})r^{-\alpha} + \frac{1}{r} \\
								&\ge \kappa\Bigl( \frac{680}{655}  + \cos\theta_{m} - \frac{2}{1000}\Bigr) > 0.
\end{align*}
Therefore, $\Re h_{m}(r\e^{\I\theta_{m}}) \le \Re h_{m}(r_{1}\e^{\I\theta_{m}}) + 2/1000^{2}$. Comparing $h_{m}(r_{1}\e^{\I\theta_{m}})$ to $h_{m}(w_{m})$, using the notations and estimates from Lemma \ref{lem: steepest path}, we get
\begin{align*}
	\Re h_{m}(r\e^{\I\theta_{m}}) &- \Re h_{m}(w_{m}) \\
	&\le \frac{\abs{h_{m}''(w_{m})}}{2}\delta^{2}\Bigl(\cos(\vphi_{m} + 2\theta_{m})\bigl((1+\Re \eta_{m}(r_{1}\e^{\I\theta_{m}})\bigr)  + \abs[1]{\Im \eta_{m}(r_{1}\e^{\I\theta_{m}})}\Bigr) + \frac{2}{1000^{2}} \\
	&\le \frac{9/10}{2}\cdot\frac{27^{2}}{680^{2}}\bigl((3/4)\cos(3\pi/4 - 8/992) + 1/4\bigr) + \frac{2}{1000^{2}} \le -\frac{1}{10000}. 
\end{align*}
We can conclude that 
\[
	\int_{\Upsilon_{3}}l_{\alpha}(\mu m^{\frac{1}{1-\alpha}}w)\e^{mh_{m}(w)}\frac{\dif w}{w} \lesssim \frac{\e^{m\Re h_{m}(w_{m})}}{\e^{m/10000}}, 
\]
which is negligible compared to \eqref{eq: contribution saddle point}.

We now let $\Upsilon_{4}$ be the arc of the circle $\abs{w-w_{m}}=\delta$ which connects $r_{1}\e^{\I\theta_{m}}$ to the initial point $c$ of $\Gamma$. Similarly, let $\Upsilon_{5}$ be the arc of the circle which connects the end point $d$ of $\Gamma$ to the point $r_{2}\e^{\I\theta_{m}}$, where $r_{2}\coloneqq \abs{w_{m}}+\delta$. Since these arcs lie in the sectors $\arg(w-w_{m}) \in (-\pi/8, \pi/8)$ and $(7\pi/8, 9\pi/8)$ respectively, the same estimate as before holds: 
\[
	\int_{\Upsilon_{4}\cup\Upsilon_{5}}l_{\alpha}(\mu m^{\frac{1}{1-\alpha}}w)\e^{mh_{m}(w)}\frac{\dif w}{w} \lesssim \frac{\e^{m\Re h_{m}(w_{m})}}{\e^{m/10000}}.
\] 

Next, we set $\Upsilon_{6}\coloneqq \{r\e^{\I\theta_{m}}: r_{2}\le r\le W\}$, with $W$ as in \eqref{eq: def eps and W}. This line segment is treated similarly as the line $\Upsilon_{3}$. We now use that the function $r\mapsto \Re f(r\e^{\I\theta_{m}})$ is \emph{non-increasing} in the range $r_{2}\le r \le W$, as apparent from a similar calculation. If $r\ge r_{2} \ge \frac{678+27}{680}\kappa^{-1}$, then
\[
	\dpd{}{r}\Re f(r\e^{\I\theta_{m}}) \le \kappa\Bigl(\frac{680}{705} + \cos\theta_{m} + \frac{1}{1000}\Bigr) < 0,
\]
since $\abs{\theta_{m}-\pi} \le 2/1000$, so $\abs{\cos\theta_{m} + 1} \le 2/1000$. Similarly as for $\Upsilon_{3}$, we conclude that 
\[
	\int_{\Upsilon_{6}}l_{\alpha}(\mu m^{\frac{1}{1-\alpha}}w)\e^{mh_{m}(w)}\frac{\dif w}{w} \lesssim \kappa W \frac{\e^{m\Re h_{m}(w_{m})}}{\e^{m/10000}}. 
\]

Finally we set $\Upsilon_{7}\coloneqq\{r\e^{\I\theta_{m}}: r\ge W\}$. Using estimate \eqref{eq: asymptotic g_{m}}, we get for $r\ge W$: 
\begin{align*}
	\Re h(r\e^{\I\theta_{m}}) 	&= \kappa r \cos\theta_{m} - \frac{1}{1-\alpha}r^{1-\alpha}\cos((1-\alpha)\theta_{m}) + \log r + \log 2 + \Re g_{m}(r\e^{\I\theta_{m}}) \\
						&= r\biggl(\kappa\cos\theta_{m} - \frac{\cos((1-\alpha)\theta_{m})}{1-\alpha}r^{-\alpha} + \frac{\log r+\log 2}{r} + O\bigl(m^{-\frac{\alpha}{1-\alpha}}r^{-\alpha}\bigr)\biggr) 
						\le -\frac{\kappa}{2}r.
\end{align*}
Hence, 
\[
	\int_{\Upsilon_{7}}l_{\alpha}(\mu m^{\frac{1}{1-\alpha}}w)\e^{mh_{m}(w)}\frac{\dif w}{w} \lesssim \int_{W}^{\infty}\e^{-m\kappa r/2}\dif r = \frac{2}{m\kappa}\e^{-m\kappa W/2}. 
\]	
This is of lower order than the contribution from the integral over $\Gamma$, since $\e^{-\kappa W/2} < c$, by the definitions of $W$ \eqref{eq: def eps and W} and $c$ \eqref{eq: def c}, and the fact that $\kappa\in I$ \eqref{eq: range kappa}. We may thus conclude that 
\begin{align*}
	\Im \int_{0}^{\I\infty}l_{\alpha}(\mu m^{\frac{1}{1-\alpha}}w)\e^{mh_{m}(w)}\frac{\dif w}{w} &= \Im\int_{\Gamma\cup \bigcup_{i}\Upsilon_{i}}l_{\alpha}(\mu m^{\frac{1}{1-\alpha}}w)\e^{mh_{m}(w)}\frac{\dif w}{w}\\
	&\gtrsim \Im \int_{\Gamma}l_{\alpha}(\mu m^{\frac{1}{1-\alpha}}w)\e^{mh_{m}(w)}\frac{\dif w}{w} \gtrsim \frac{c^{m}}{\sqrt{m}},
\end{align*}
as claimed earlier.	

\mbox{}

%\begin{acknowledgements} 
The authors thank Prof.\ Jasson Vindas for supporting the research project which led to this paper, and initiating discussion between the authors, resulting in a fruitful collaboration.
%\end{acknowledgements}

\end{document}